\documentclass[aos]{imsart}

\RequirePackage{amsthm,amsmath,amsfonts,amssymb}
\RequirePackage[authoryear]{natbib}
\RequirePackage[colorlinks,citecolor=blue,urlcolor=blue]{hyperref}

\startlocaldefs
\newtheorem{lem}{Lemma}
\newtheorem{theo}{Theorem}
\newtheorem{rem}{Remark}

\endlocaldefs

\begin{document}

\begin{frontmatter}
\title{Asymptotic Equivalence of Locally Stationary Processes and Bivariate Gaussian White Noise}
\runtitle{Asymptotic Equivalence of Locally Stationary Processes}

\begin{aug}
\author[A]{\fnms{Cristina}~\snm{Butucea}\ead[label=e1]{cristina.butucea@ensae.fr}},
\author[B]{\fnms{Alexander}~\snm{Meister}\ead[label=e2]{alexander.meister@uni-rostock.de}}
\and
\author[C]{\fnms{Angelika}~\snm{Rohde}\ead[label=e3]{angelika.rohde@stochastik.uni-freiburg.de}}
\address[A]{CREST, ENSAE, Institut Polytechnique de Paris, 5, avenue Henry Le Chatelier, Palaiseau, 91120, France\printead[presep={,\ }]{e1}}

\address[B]{Institut f\"ur Mathematik, Universit\"at Rostock, Ulmenstra{\ss}e 69, 18057 Rostock, Germany\printead[presep={,\ }]{e2}}

\address[C]{Mathematisches Institut, Albert-Ludwigs-Universit\"at Freiburg, Ernst-Zermelo-Stra{\ss}e 1, 79104 Freiburg, Germany\printead[presep={,\ }]{e3}}
\end{aug}

\begin{abstract}
We consider a general class of statistical experiments, in which an $n$-dimensional centered Gaussian random variable is observed and its covariance matrix is the parameter of interest. The covariance matrix is assumed to be well-approximable in a linear space of lower dimension $K_n$ with eigenvalues uniformly bounded away from zero and infinity. We prove asymptotic equivalence of this experiment and a class of $K_n$-dimensional Gaussian models with informative expectation in Le Cam's sense when $n$ tends to infinity and $K_n$ is allowed to increase moderately in $n$ at a polynomial rate. For this purpose we derive a new localization technique for non-i.i.d. data and a novel high-dimensional Central Limit Law in total variation distance. These results are key ingredients to show asymptotic equivalence between the experiments of locally stationary Gaussian time series and a bivariate Wiener process with the log spectral density as its drift. On this way a novel class of matrices is introduced which generalizes circulant Toeplitz matrices traditionally used for strictly stationary time series.
\end{abstract}



\end{frontmatter}


\section{Introduction}

Suppose that an $n$-dimensional centered Gaussian random vector $X = (X_1,\ldots,X_n)$ is observed, where the $n\times n$-covariance matrix $\theta$ is unknown to the statistician. Thus the corresponding statistical experiment equals
$$ {\cal A}_n \, := \, \big(\mathbb{R}^n, \mathfrak{B}(\mathbb{R}^n), \{{\cal N}(0,\theta) : \theta\in\Theta\}\big)\,, $$
with the $n$-dimensional Borel $\sigma$-field $\mathfrak{B}(\mathbb{R}^n)$. The parameter space $\Theta$ is assumed to consist of symmetric and positive definite $n\times n$-matrices with eigenvalues uniformly bounded away from zero and infinity and  well-approximable in a linear space of lower dimension $K_n$.  Our goal of this paper is two-fold:
\begin{itemize}
    \item[(i)] We aim at  developing  a widely applicable theory for deriving asymptotic equivalence of $(\mathcal{A}_n)$ and a sequence of accompanying Gaussian mean models in the sense of \cite{LeCam1964}, laying the foundation to move over to central nonparametric statistical models in the presence of intricate dependency structures.
\item[(ii)] Specifying $\Theta$ to covariance matrices of locally stationary processes, we finalize our general strategy in (i) with  proving asymptotic equivalence for time-varying  spectral density estimation of  locally stationary processes and bivariate Gaussian white noise. \cite{GolNussZhou2010} showed asymptotic equivalence of stationary processes and univariate Gaussian white noise. In contrast, our work is providing a substantial generalization of the strictly stationary case with a new and significantly different methodology.
\end{itemize}
A family of probability distributions $\mathcal{P}=\{\mathbb{P}_{\theta}\arrowvert\, \theta\in \Theta\}$ on some measurable space $(\mathcal{X},\mathcal{A})$ is called statistical model or statistical experiment. For two experiments $\mathcal{P}=\{\mathbb{P}_{\theta}\arrowvert\, \theta\in \Theta\}$ and $\mathcal{Q}=\{\mathbb{Q}_{\theta}\arrowvert \theta\in\Theta\}$, indexed by the same parameter set but defined  on potentially different measurable spaces $(\mathcal{X},\mathcal{A})$ and $(\mathcal{Y},\mathcal{B})$, respectively, \cite{LeCam1964} has introduced the concept of deficiency $\delta(\mathcal{P},\mathcal{Q})$ which quantifies statistical approximability of $\mathcal{Q}$ by $\mathcal{P}$. The most common way of bounding the deficiency $\delta(\mathcal{P},\mathcal{Q})$ is to propose a Markov kernel $K$ which transfers distributions from $(\mathcal{X},\mathcal{A})$ to $(\mathcal{Y},\mathcal{B})$  such that  $\sup_{\theta\in \Theta}\Arrowvert K\mathbb{P}_{\theta}-\mathbb{Q}_{\theta}\Arrowvert_{TV}$ is small. Two sequences $(\mathcal{P}_n)$ and $(\mathcal{Q}_n)$ of statistical experiments are called asymptotically equivalent if $\Delta(\mathcal{P}_n,\mathcal{Q}_n)=\max\big(\delta(\mathcal{P}_n,\mathcal{Q}_n),\delta(\mathcal{Q}_n,\mathcal{P}_n)\big)\rightarrow 0$. In a nutshell, asymptotic equivalence means that any asymptotic solution to a statistical decision problem in one experiment automatically reveals a corresponding approximate solution in the other one and vice versa. Moreover, if these results are constructive, there is an explicit recipe for this correspondence. Due to its overriding statistical importance, asymptotic equivalence of nonparametric experiments has constantly received  a lot of attention in the literature, see \cite{Nussbaum96} for density estimation, \cite{BrownLow1996} and \cite{Rohde+2004+235+243} for nonparametric regression, \cite{JaehnischNussbaum2003} for independent but not identically distributed data, \cite{GramaNeumann2006} for autoregression, \cite{DalalyanReiss2006} and \cite{DalalyanReiss2007} for diffusion models, \cite{Carter2007} for regression with unknown conditional variance, \cite{Reiss2008} for multivariate regression, \cite{GolNussZhou2010} for Gaussian strictly stationary time series, \cite{Meister2011} for functional linear regression, \cite{Reiss2011} for volatility estimation based on noisy discrete-time observations of a continuous martingale, \cite{SchmidtHieber2014} for regression under fractional noise, \cite{ButGutaNussbaum2018} for quantum statistical models.  As a consequence of the concept's nature, constructive proofs are inherently case-specific. Hence, it has become time to think of a unifying framework which is applicable to a variety of models at once.

Towards this goal, a general Gaussianization scheme for $\mathcal{A}_n$ is developed in Section~\ref{GGS}, revealing  asymptotic equivalence of $(\mathcal{A}_n)$ and a class of $K_n$-dimensional Gaussian models with informative expectation. At the core of this scheme are a new localization technique for potentially non-i.i.d.~data as well as  a new quantitative high-dimensional central limit theorem for Gaussian quadratic forms in total variation distance  with explicit dependence on the dimension.

Based on this Gaussianization scheme, the problem of asymptotic equivalence for Gaussian locally stationary processes and a bivariate Gaussian white noise model is addressed in Section~\ref{sec:AE_LSP_BGW}. Locally stationary processes in the sense of \cite{Dahlhaus1996,Dahlhaus1997,10.1214/aos/1015957480} provide a theoretical framework for rigorous local asymptotic statements beyond stationarity. Since then and until now, they have been of great interest in the statistical literature, see \cite{10.1214/19-AOS1895}, \cite{10.1214/22-AOS2205}, \cite{10.1214/24-AOS2358}, \cite{10.1214/23-AOS2288}, \cite{10.1214/18-AOS1743}, \cite{T-G_Nussbaum2012}, \cite{10.1214/12-AOS1043}, \cite{10.1214/20-AOS1954} to mention just a few recent ones; see also \cite{10.1214/18-AOS1739} for the related concept with Markov chains. Locally stationary processes  are given as triangular array of processes with time-varying spectral representation whose Wigner-Ville spectra converge under suitable smoothness assumptions towards a time-varying spectral density in quadratic mean. We establish asymptotic equivalence between the experiments of locally stationary Gaussian time series and a bivariate Wiener process with the time-varying log spectral density as its drift. Therein a novel class of matrices is introduced which generalizes circulant Toeplitz matrices traditionally used for strictly stationary time series.

Beyond local stationarity, our general results of Section~\ref{GGS} could find applications to other Gaussian process models as, for instance,  described in \cite{RasmussenWilliams} and \cite{vdVaart_vZanten}.

\newpage

\section{General Gaussianization Scheme} \label{GGS}

\subsection{Pre-Smoothing} \label{PreSmooth}

 At this stage we impose that all eigenvalues of the parameter $\theta\in \Theta$ are bounded from below by $\rho$ and from above by $1/\rho$ for some fixed constant $\rho>0$ -- uniformly with respect to $\theta \in \Theta$. This also implies that $\rho \, \leq \, 1$. The first step is to approximate the original experiment ${\cal A}_n$ by another Gaussian experiment ${\cal B}_n$ where the covariance matrix $C_\theta$ stems from a specific known $K_n$-dimensional linear subspace $L_n$ of the space of all symmetric $n\times n$-matrices. Therein $K_n$ is assumed to be much smaller than the dimension $N:=n(n+1)/2$ of the original space of the symmetric $n\times n$-matrices. Let $M_1,\ldots,M_{K_n}$ denote any orthonormal basis of $L_n$ with respect to the Frobenius inner product. The covariance matrix of the pre-smoothed Gaussian experiment equals the projected matrix
$$ C_\theta \, = \, \sum_{k=1}^{K_n} \langle \theta,M_k\rangle_F \cdot M_k\,, $$
where we write $\alpha_{\theta,k} := \langle \theta,M_k\rangle_F$ and $\alpha_\theta := (\alpha_{\theta,1},\ldots,\alpha_{\theta,K_n})^\top$; and $\langle\cdot,\cdot\rangle_F$ stands for the Frobenius inner product. This general approach is applicable to the covariance matrix of locally stationary time series by specific choice of the basis matrices $M_1,\ldots,M_{K_n}$ (Section \ref{LocStat}). We apply the equation (A.4) in \cite{Reiss2011} to establish asymptotic equivalence of the original experiment ${\cal A}_n$ and
$$ {\cal B}_n \, := \, \big(\mathbb{R}^n, \mathfrak{B}(\mathbb{R}^n), \{{\cal N}(0,C_\theta) : \theta\in\Theta\}\big)\,, $$
whenever
\begin{equation} \label{eq:cond1} \lim_{n\to\infty} \, \sup_\theta \, \big\|\theta^{-1/2} (\theta - C_\theta) \theta^{-1/2}\big\|_F^2 \, = \, 0\,, \end{equation}
where $\|\cdot\|_F$ denotes the Frobenius norm. We deduce that
\begin{align} \nonumber
\big\| \theta^{-1/2}  (\theta - C_\theta) \theta^{-1/2}\big\|_F^2  & \, \geq \, \big\|I_n - \theta^{-1/2} C_\theta \theta^{-1/2}\big\|_{sp}^2  
\\ & \label{eq:cond1.1}
\, \geq \, \sup_{\|w\|=1} \, \big|1 - (w^\top C_\theta w) / (w^\top \theta w)\big|^2,
\end{align}
where $\|\cdot\|_{sp}$ and $I_n$ denote the spectral norm and the $n\times n$-identity matrix, respectively. Therefore, validity of (\ref{eq:cond1}) guarantees that, for any constant $c \in (0,1)$, all eigenvalues of $C_\theta$ have the lower bound $c\rho$ for all $\theta\in\Theta$ when $n$ is sufficiently large. This also implies positive definiteness of $C_\theta$. Analogously we derive that, for any constant $c \in (0,1)$, all eigenvalues of $C_\theta$ have the upper bound $1/(c\rho)$ for all $\theta\in\Theta$ when $n$ is sufficiently large.

\subsection{A Novel Approach to Localization} \label{GLS}

We suggest a novel general localization scheme which avoids splitting the sample as that technique seems awkward for experiments with non i.i.d. data. Consider some general experiment ${\cal X}_n$ in which some random variable $X$ is observed with ${\cal L}(X) = \mathbb{P}_\theta$, $\theta \in\Theta$, with $\Theta$ as before. In particular, if we take $\mathbb{P}_\theta$ to be the Gaussian vector $\mathcal{N}(0,\theta)$ we get the initial experiment $\mathcal{A}_n$. Let $\hat{\alpha}$ be any estimator of the vector $\alpha_\theta$ from the previous section based on $X$. Obviously, ${\cal X}_n$ is equivalent to the experiment ${\cal X}'_n$ which describes the observation of $(X,\hat{\alpha}(X)+\eta)$ where $X$ and $\eta$ are independent and ${\cal L}(\eta) = {\cal N}(0,\beta_n^2 I_{K_n})$ where $\beta_n$ is some scaling factor to be selected. Thus the estimator $\hat{\alpha}$ is artificially and intentionally contaminated by the Gaussian random variable $\eta$. We derive an upper bound on the Le Cam distance between ${\cal X}'_n$ and the experiment ${\cal X}''_n$ in which $(X,\alpha_\theta + \eta)$ is observed.
\begin{align*}
\frac12 \, & \sup_\theta \sup_{\|\Phi\|_\infty\leq 1} \big|\mathbb{E}_\theta \Phi(X,\hat{\alpha}(X)+\eta) - \mathbb{E}_\theta \Phi(X,\alpha_\theta+\eta)\big| \\
& \, \leq \, \frac12 \, \sup_\theta \, \mathbb{E}_\theta  \sup_{\|\Phi\|_\infty\leq 1} \big|\mathbb{E}_\theta \big\{\Phi(X,\hat{\alpha}(X)+\eta) \mid X\big\} - \mathbb{E}_\theta \big\{\Phi(X,\alpha_\theta+\eta) \mid X\big\}\big| \\
& \, = \, \sup_\theta \mathbb{E}_\theta \, \mbox{TV}\big\{{\cal L}\big(\hat{\alpha}(X)+\eta\mid X\big) , {\cal L}(\alpha_\theta + \eta) \mid X\big\} \\
& \, \leq \, \big\{1 - \exp\big[ - \sup_\theta \, \mathbb{E}_\theta \, {\cal K}\big({\cal N}(\hat{\alpha}(X),\beta_n^2 I_{K_n}) \, , \, {\cal N}(\alpha_\theta,\beta_n^2 I_{K_n}) \mid X\big)\big]\big\}^{1/2} \\
& \, = \, \big\{1 - \exp\big[- \beta_n^{-2} \cdot \sup_\theta \, \mathbb{E}_\theta \|\hat{\alpha} - \alpha_\theta\|^2\big]\big\}^{1/2}\,,
\end{align*}
where $\mbox{TV}(P,Q\mid X)$ and ${\cal K}(P,Q\mid X)$ denote the conditional total variation distance and the conditional Kullback-Leibler divergence between some random probability measures $P$ and $Q$ given $X$, respectively. The Bretagnolle-Huber inequality as well as Jensen's inequality have been used. Thus,
\begin{equation} \label{eq:cond2}
\sup_\theta \, \mathbb{E}_\theta \|\hat{\alpha} - \alpha_\theta\|^2 \, = \, o(\beta_n^2)\,, \end{equation}
implies asymptotic equivalence of ${\cal X}_n'$ and ${\cal X}_n''$. Finally we introduce the experiment ${\cal X}_n'''$ in which $(X,\alpha_\theta+\tilde{\eta})$ is observed where $X$ and $\tilde{\eta}$ are independent and $\tilde{\eta}$ has the truncated density
$$ f_{\tilde{\eta}} \, := \, f_\eta \cdot 1_{[0,\gamma_n]}(\|\cdot\|) \, \bigg/ \, \int f_\eta(x) \cdot 1_{[0,\gamma_n]}(\|x\|) \, dx\,, $$
for some sequence $(\gamma_n)_n$ and the ${\cal N}(0,\beta_n^2 I_{K_n})$-density $f_\eta$. The Le Cam distance between ${\cal X}_n''$ and ${\cal X}_n'''$ obeys the upper bound
\begin{align*}
\frac12 \,  \sup_\theta \, \int  \big|f_\eta(x-\alpha_\theta) - f_{\tilde{\eta}}(x-\alpha_\theta)\big| \, dx &
\, \leq \, \int \big|f_\eta(x) - f_\eta(x)\cdot 1_{[0,\gamma_n]}(\|x\|)\big|\, dx \\
& \, = \, \mathbb{P}\big[\|\eta\|^2 > \gamma_n^2\big]  
\, \leq \, K_n \cdot \beta_n^2 / \gamma_n^2\,, \end{align*}
where 
we have used Markov's inequality in the last step. Hence, imposing (\ref{eq:cond2}) and
\begin{equation} \label{eq:cond3}
\lim_{n\to\infty} \, K_n \cdot \beta_n^2 / \gamma_n^2 \, = \, 0\,,
\end{equation}
we have asymptotic equivalence between ${\cal X}_n$ and ${\cal X}_n'''$. Thanks to the independence of the observed components in ${\cal X}_n'''$, one can apply $\alpha := \alpha_\theta + \tilde{\eta}$ in order to localize the experiment ${\cal X}_n$, i.e. we can treat and transform the experiment ${\cal X}_n$ under the information that, almost surely, the (Euclidean) distance between some additional independent observation $\alpha_\theta + \tilde{\eta}$ and the true $\alpha_\theta$ is bounded from above by $\gamma_n$.

\subsection{Dimension Reduction in the Localized and Smoothed Experiment} \label{DR}

In our concrete experiment ${\cal B}_n$ from Section \ref{PreSmooth}, we specify the symmetric $n\times n$-matrix
\begin{equation} \label{eq:C} C := \sum_{k=1}^{K_n} \big(\alpha_{\theta,k} + \tilde{\eta}_k\big) \cdot M_k\,, \end{equation}
as the localized version of the unknown matrix $C_\theta$. Thus, by ${\cal C}_n$, we denote the experiment, in which the independent random localization matrix $C$ (or -- equivalently -- the vector $\alpha := \{\alpha_{\theta,k} + \tilde{\eta}_k\}_{k=1,\ldots,K_n}$ by the arguments from Section \ref{GLS}) is available in addition to the observation from experiment ${\cal B}_n$. It is asymptotically equivalent to ${\cal B}_n$ whenever (\ref{eq:cond2}) and (\ref{eq:cond3}) hold true. As a pilot estimator in the experiment ${\cal B}_n$ we employ $$ \hat{\alpha}(x) \, := \, \big\{\langle x x^\top, M_k\rangle_F\big\}_k\,, $$
for all $x\in \mathbb{R}^n$. We deduce that
\begin{align*}
\mathbb{E}_\theta \big\|\hat{\alpha} - \alpha_\theta\big\|^2 & \, = \, 2 \sum_{k=1}^{K_n} \big\|C_\theta^{1/2} M_k C_\theta^{1/2}\big\|_F^2 \, \leq \, 2 K_n \|C_\theta\|_{sp}^2 \, \leq \, 4 K_n \rho^{-2}\,,
\end{align*}
uniformly for all $\theta\in \Theta$ when $n$ is sufficiently large. We have used that $\|AB\|_F \leq \|A\|_{sp} \|B\|_F$ and $\|AB\|_F \leq \|A\|_F \|B^\top\|_{sp}$ for all $n\times n$-matrices $A$ and $B$.  Therefore (\ref{eq:cond2}) and (\ref{eq:cond3}) are satisfied for appropriate selection of $\beta_n$ if \begin{equation} \label{eq:cond.pilot1}
\lim_{n\to\infty} \, K_n^2\, \gamma_n^{-2} \, = \, 0\,.
\end{equation}
Then, for $\Delta_\theta := C - C_\theta$, it holds that
\begin{align} \nonumber
\sup_\theta \, \|\Delta_\theta\|_F^2 & \, \leq \, \gamma_n^2\,, \\ \label{eq:gamman}
\sup_\theta \, \|\Delta_\theta\|_{sp}^2 & \, \leq \, \gamma_n^2 \cdot \sum_{k=1}^{K_n} \|M_k\|_{sp}^2\,.
\end{align}
Now we impose the crucial condition
\begin{equation} \label{eq:condMk}
\sup_n \, \sup_{k=1,\ldots,K_n}  \, n \cdot \|M_k\|_{sp}^2 \, < \, \infty\,.
\end{equation}
Then, under the constraint
\begin{equation} \label{eq:cond4}
\lim_{n\to\infty}\, \gamma_n^{2} \, K_n / n \, = \, 0\,,
\end{equation}
all eigenvalues of $C$ are bounded from below by $c\rho$ and from above by $1/(c\rho)$ for any constant $c\in (0,1)$ when $n$ is sufficiently large (uniformly with respect to $\theta$). Furthermore the matrix $C$ is invertible and all of its eigenvalues also admit the lower bound $c \rho$ for any constant $c\in (0,1)$ for $n$ large enough (again uniformly in $\theta$); moreover,
\begin{align*}
C_\theta^{-1} & \, = \, C^{-1} + C_\theta^{-1} - C^{-1} \, = \, C^{-1} + C^{-1}(C C_\theta^{-1} - I_n) \, = \, C^{-1} + C^{-1} \Delta_\theta C_\theta^{-1} \\
& \, = \, C^{-1} + C^{-1} \Delta_\theta \big(C^{-1} + C^{-1} \Delta_\theta C_\theta^{-1}\big) \, = \, B_\theta \, + \, C^{-1} \Delta_\theta C^{-1} \Delta_\theta C_\theta^{-1}\,,
\end{align*}
where
$$ B_\theta \, := \, C^{-1} + C^{-1} \Delta_\theta C^{-1} \, = \, \big(I_n + C^{-1} \Delta_\theta\big) C^{-1}\,. $$
We realize that the matrix $B_{\theta}$ is symmetric and invertible whenever
\begin{equation*}
\sup_\theta \, \big\|C^{-1} \Delta_\theta\big\|_{sp} \, < \, 1\,,
\end{equation*}
which follows from (\ref{eq:cond4}). Then the inverse of $B_{\theta}$ may be expanded by the Neumann series; concretely,
\begin{align*}
B_{\theta}^{-1} & \, = \, C \sum_{k=0}^\infty \big(-C^{-1} \Delta_\theta\big)^{k} \, = \, C_\theta \, + \, C \sum_{k=2}^\infty \big(-C^{-1} \Delta_\theta\big)^{k}\,.
\end{align*}
We deduce that
\begin{align*}
\big\| & C_\theta^{-1/2} \big(B_{\theta}^{-1} - C_\theta\big) C_\theta^{-1/2} \big\|_F  \, = \, \Big\| \sum_{k=2}^\infty (-1)^k C_\theta^{-1/2} \Delta_\theta  (C^{-1} \Delta_\theta)^{k-1} C_\theta^{-1/2}\Big\|_F \\
& \, \leq \, \sum_{k=1}^\infty \|C_\theta^{-1/2}\|_{sp}^2 \, \|C^{-1}\|_{sp}^k \, \|\Delta_\theta\|_{sp}^k \, \|\Delta_\theta\|_F \, \leq \, \sum_{k=1}^\infty  \{\gamma_n /(c \rho)\}^{k+1} \, \Big(\sum_{k=1}^{K_n} \|M_k\|_{sp}^2\Big)^{k/2}\,,
\end{align*}
for any $c\in (0,1)$ when $n$ is sufficiently large. This upper bound converges to $0$ under (\ref{eq:condMk}) and
\begin{equation} \label{eq:cond4'}
\lim_{n\to\infty} \, \gamma_n^4 \, K_n / n \, = \, 0\,,
\end{equation}
which represents a stronger condition than (\ref{eq:cond4}) as we grant (\ref{eq:cond.pilot1}). Using (A.4) from \cite{Reiss2011} again, we show that  (\ref{eq:condMk}) and (\ref{eq:cond4'}) imply asymptotic equivalence of the experiments ${\cal C}_n$ and ${\cal D}_n$, in which one observes $C$ and -- conditionally on $C$ -- a ${\cal N}(0,B_\theta^{-1})$-distributed random variable. Note that, using the arguments of (\ref{eq:cond1.1}), the constraint (\ref{eq:cond4'}) also provides that
$$ \lim_{n\to\infty} \, \sup_\theta \, \sup_{\|v\|=1} \, \big|1 \, - \, v^\top B_{\theta}^{-1} v \, / \, v^\top C_\theta v\big|^2 \, = \, 0\,, $$
such that, for any constant $c\in (0,1)$, the smallest eigenvalue of $B_{\theta}^{-1}$ is larger or equal to $c\rho$ for all $\theta$ when $n$ is sufficiently large. This guarantees positive definiteness of $B_{\theta}^{-1}$ for these $n$.

For any $x = (x_1,\ldots,x_n)^\top$, consider the conditional likelihood function of the experiment ${\cal D}_n$ given $C$ with respect to the $n$-dimensional Lebesgue-Borel measure
\begin{align*} f_n&(\theta;x\mid C) \, = \, (2\pi)^{-n/2} \, (\det B_{\theta})^{1/2} \, \exp\big\{-x^\top B_{\theta} x / 2\big\} \\ &
\, = \,  (2\pi)^{-n/2} \, (\det B_{\theta})^{1/2} \, \exp\big\{-x^\top C^{-1} \Delta_\theta C^{-1} x / 2\Big\} \cdot \exp\big\{ - x^\top C^{-1} x/2\big\} \\ &
\, = \,  (2\pi)^{-n/2} \, (\det B_{\theta})^{1/2} \, \exp\Big\{-\frac12\sum_{k=1}^{K_n} \langle\Delta_\theta , M_k\rangle_F \cdot x^\top C^{-1} M_k C^{-1} x\Big\} \\ & \hspace{9cm} \cdot \exp\big\{ - x^\top C^{-1} x/2\big\}\,.
\end{align*}
By Fisher-Neyman factorization, we have shown that $(C,T)$ with
$$ T(x) \, := \, \big\{x^\top C^{-1} M_k C^{-1} x\big\}_{k=1,\ldots,K_n} $$
forms an sufficient statistic in the experiment ${\cal D}_n$. Therefore, ${\cal D}_n$ is equivalent to the experiment ${\cal E}_n$ in which only the statistic $(C,T)$ is observed.

Moreover, the experiment ${\cal E}_n$ is also asymptotically equivalent to the experiment ${\cal F}_n$ which describes the observation of the statistic $(C,T)$ when $T$ is considered under the basic probability measure ${\cal N}(0,C_\theta)$ -- instead of ${\cal N}(0,B_{\theta}^{-1})$ as in ${\cal E}_n$. This follows directly from the asymptotic equi\-valence of ${\cal C}_n$ and ${\cal D}_n$ (without any transforming Markov kernels) that has already been shown. Note that the total variation distance between two probability measures $P$ and $Q$ dominates that between the corresponding image measures under the same mapping.

\subsection{High-dimensional Central Limit Theorem in Total Variation Distance}

The observed component $T$ in the experiment ${\cal F}_n$ may be represented by
\begin{equation}\label{Tvec}
    \varepsilon^\top C_\theta^{1/2} C^{-1} M_k C^{-1} C_\theta^{1/2} \varepsilon\,, \quad k=1,\ldots,K_n\,,
\end{equation}
where $\varepsilon = (\varepsilon_1,\ldots,\varepsilon_n)^\top$ is a random vector with independent ${\cal N}(0,1)$-distributed components that is also independent of $C$. Let us consider the conditional characteristic function $\Psi_n$ of $T$ given $C$.
\begin{align} \nonumber
\Psi_n(t_1,\ldots,t_{K_n}) & \, = \,  \mathbb{E}^C \exp\big( i \varepsilon^\top C_\theta^{1/2} C^{-1} M^{[t]} C^{-1} C_\theta^{1/2} \varepsilon\big) \, = \, \mathbb{E}^C \exp\Big( i \sum_{j=1}^n \lambda_{j,\theta,t} \, \varepsilon_j^2\Big)\\ \label{eq:Gauss.1} &
 \, = \, \prod_{j=1}^n \big(1 -2i\lambda_{j,\theta,t}\big)^{-1/2} \, = \, \big\{\det\big(I_n -2i D_\theta^{[t]}\big)\big\}^{-1/2}\,,
\end{align}
for all $t=(t_1,\ldots,t_{K_n})^\top \in \mathbb{R}^{K_n}$ where $M^{[t]} := \sum_{k=1}^{K_n} t_k M_k$; the $\lambda_{j,\theta,t}$ denote the eigenvalues of the symmetric matrix $D_\theta^{[t]} := C_\theta^{1/2} C^{-1} M^{[t]} C^{-1} C_\theta^{1/2}$; and, here, we write $\mathbb{E}^C$ for the conditional expectation given $C$. We have used the invariance of the distribution of $\varepsilon$ with respect to the action of orthogonal matrices; the fact that each $\varepsilon_j^2$ has the $\chi^2(1)$-distribution; and the coincidence of the determinants of similar matrices. It follows from Fredholm's representations of the determinant (see e.g. \cite{Fredholm1903}, p. 384; \cite{Bornemann2010}) that
\begin{align} \nonumber
\Psi_n(t_1,\ldots,t_{K_n}) & \, = \, \exp\Big\{ \frac12 \sum_{\ell=1}^\infty \frac1\ell (2i)^\ell \cdot \mbox{tr}\big[ \big(D_\theta^{[t]}\big)^\ell\big]\Big\} \\ \label{eq:CL1}
& \, = \, \exp\Big\{ i \sum_{k=1}^{K_n} t_k  d_{\theta,k} - \frac12 t^\top \Gamma_\theta t\Big\} \cdot \exp\Big\{ \frac12 \sum_{\ell=3}^\infty \frac1\ell (2i)^\ell \cdot \mbox{tr}\big[ \big(D_\theta^{[t]}\big)^\ell\big]\Big\}\,,
\end{align}
if $\|D_\theta^{[t]}\|_{sp} < 1/2$ where we write tr for the trace of a matrix; and
\begin{align*}
d_{\theta,k} & \, := \, \mbox{tr}\big(C_\theta C^{-1} M_k C^{-1}\big)\,, \\
\Gamma_\theta & \, := \, 2 \, \big\{\mbox{tr}\big(C_\theta C^{-1} M_k C^{-1} C_\theta C^{-1} M_{k'} C^{-1}\big)\big\}_{k,k'} \\
& \, = \, 2 \, \big\{\big\langle C_\theta^{1/2} C^{-1} M_k C^{-1} C_\theta^{1/2} , C_\theta^{1/2} C^{-1} M_{k'} C^{-1} C_\theta^{1/2} \big\rangle_F\big\}_{k,k'}\,.
\end{align*}
Note that the first factor in (\ref{eq:CL1}) represents the characteristic function of a ${\cal N}\big(\{d_{\theta,k}\}_k , \Gamma_\theta\big)$-distributed random variable at $t$. The total variation distance between ${\cal L}_\theta(T|C)$ and ${\cal N}(\{d_{\theta,k}\}_k,\Gamma_\theta)$ equals
$$ \mbox{TV}\big({\cal L}_\theta(T|C)\, , \, {\cal N}(\{d_{\theta,k}\}_k,\Gamma_\theta)\big) \, = \, \mbox{TV}\big({\cal L}_\theta(T^*|C)\, , \,{\cal N}(0,I_{K_n})\big)\,, $$
where $T^* := \Gamma_\theta^{-1/2} (T - \{d_{\theta,k}\}_k)$ as the total variation distance is invariant with respect to shifts and the action of invertible matrices. The positive definiteness of $\Gamma_\theta$ is guaranteed by the following lemma, where we define
\begin{align} \nonumber
\Gamma & \, := \, 2\, \big\{\mbox{tr}(M_k C^{-1} M_{k'} C^{-1})\big\}_{k,k'}\,, \\ \label{eq:tildeGamma}
\tilde{\Gamma}_\theta & \, := \, 2\, \big\{\mbox{tr}(M_k C_\theta^{-1} M_{k'} C_\theta^{-1})\big\}_{k,k'}\,.
\end{align}
The following lemma provides the exchangeability of these matrices as covariance matrices of asymptotically equivalent experiments.

\begin{lem} \label{L:Gamma_theta} Grant (\ref{eq:condMk}) and assume that
$\lim_{n\to\infty} \, K_n^2 \, \gamma_n^2 / n\, = \, 0$. 
Then
$$ \lim_{n\to\infty} \, \sup_\theta \big\{\big\|\Gamma^{-1/2} (\Gamma-\Gamma_\theta) \Gamma^{-1/2}\big\|_F \, + \, \big\|\Gamma^{-1/2} (\Gamma-\tilde{\Gamma}_\theta) \Gamma^{-1/2}\big\|_F\big\} \, = \, 0\,. $$
Moreover all eigenvalues of $\Gamma$, $\Gamma_\theta$ and $\tilde{\Gamma}_\theta$ admit the lower bound $c \rho^2$
for any fixed $c\in (0,2)$ when $n$ is sufficiently large (uniformly with respect to $\theta$).
\end{lem}

\noindent Note that the conditions of Lemma \ref{L:Gamma_theta} follow from (\ref{eq:condMk}) and (\ref{eq:cond4'}). The conditional characteristic function $\Psi_n^*$ of $T^*$ given $C$ turns out to be
$$ \Psi_n^*(t_1,\ldots,t_{K_n}) \, = \, \exp\Big\{ -i \sum_{k=1}^{K_n} \big(\Gamma_\theta^{-1/2}t\big)_k \, d_{\theta,k} \Big\} \cdot \Psi_n\big(\Gamma_\theta^{-1/2} t\big)\,, \qquad t=(t_1,\ldots,t_{K_n})^\top\,, $$
so that, by (\ref{eq:Gauss.1}),
\begin{equation} \label{eq:psin*}
\Psi_n^*(t_1,\ldots,t_{K_n}) \, = \, \exp\big(-\|t\|^2/2\big) \cdot \exp\Big\{\frac12 \sum_{\ell=3}^\infty \frac1\ell (2i)^\ell \cdot \mbox{tr}\big[\big(\tilde{D}_\theta^{[t]}\big)^\ell\big]\Big\}\,,
\end{equation}
for $\big\|\tilde{D}_\theta^{[t]}\big\|_{sp} < 1/2$ where $\tilde{D}_\theta^{[t]} := \sum_{k=1}^{K_n} t_k \cdot \tilde{D}_{\theta,k}$ and
$$ \tilde{D}_{\theta,k} \, := \, C_\theta^{1/2} C^{-1} \sum_{k'=1}^{K_n} \big(\Gamma_\theta^{-1/2}\big)_{k,k'} M_{k'} \, C^{-1} C_\theta^{1/2}\,. $$
The symmetric $n\times n$-matrices $\tilde{D}_{\theta,k}$, $k=1,\ldots,K_n$, form an orthogonal system with respect to the Frobenius inner product and satisfy $\|\tilde{D}_{\theta,k}\|_F^2=1/2$ for all $k=1,\ldots,K_n$. Moreover consider that
\begin{align} \nonumber
\big\|\tilde{D}_\theta^{[t]}\big\|_{sp} & \, \leq \, \|C_\theta^{1/2} C^{-1}\|_{sp} \cdot \|C^{-1} C_\theta^{1/2}\|_{sp} \cdot \Big\| \sum_{k=1}^{K_n} \big(\Gamma_\theta^{-1/2} t\big)_k M_k\Big\|_{sp} \\ \nonumber
& \, \leq \, \|C_\theta^{1/2} C^{-1}\|_{sp} \cdot \|C^{-1} C_\theta^{1/2}\|_{sp} \cdot\sum_{k=1}^{K_n} \big|(\Gamma_\theta^{-1/2} t)_k\big|\cdot \|M_k\|_{sp} \\ \nonumber
& \, \leq \, \|t\| \cdot \|C_\theta^{1/2} \|_{sp}^2 \cdot \|C^{-1}\|_{sp}^2  \cdot \|\Gamma_\theta^{-1/2}\|_{sp} \cdot \Big(\sum_{k=1}^{K_n} \|M_k\|_{sp}^2\Big)^{1/2}  \\ \label{eq:bound_Dt}
& \, \leq \, \|t\| \cdot \mu_n\,,
\end{align}
for $n$ sufficiently large (uniformly with respect to $\theta$), where
\begin{align*} \mu_{n} & \, \asymp \, K_n^{1/2} n^{-1/2}\,,
\end{align*}
with a constant factor (also uniform with respect to $\theta$), by (\ref{eq:condMk}) and Lemma \ref{L:Gamma_theta} so that validity of the expansion (\ref{eq:psin*}) is guaranteed whenever $\|t\| < 1/(2\mu_{n})$. Moreover (\ref{eq:cond.pilot1}) and (\ref{eq:cond4'}) also yield that
$\lim_{n\to\infty} \mu_{n} \, = \, 0$.

We define the experiment ${\cal G}_n$, in which one observes $C$ and -- conditionally on $C$ -- a ${\cal N}(\{d_{\theta,k}\}_k , \Gamma_\theta )$-distributed random vector. In the following subsections we establish asymptotic equivalence between the experiments ${\cal F}_n$ and ${\cal G}_n$ under certain conditions.

\subsubsection{Upper Bound on the Fourier Tails}

In (\ref{eq:psin*}) we can approximate $\Psi_n^*$ by a Gaussian characteristic function on a bounded domain. Still, as we are seeking to bound the total variation distance, the tails of the characteristic function $\Psi_n^*$ need to be considered as well. They are studied in the following lemma.
\begin{lem}  \label{L:FouTails}
Impose that $K_n\geq 2$. If $\mu_{n}^{-2} > 8 K_n$ then the function $\Psi_n^*$ from (\ref{eq:psin*}) is absolutely integrable over the whole of $\mathbb{R}^{K_n}$. If even $\mu_{n}^{-2}> 8K_n+16$, then it holds that
\begin{align*} \idotsint_{\|t\|\geq R} |\Psi_n^*(t_1,\ldots,t_{K_n})| & \, dt_1\cdots dt_{K_n} \\ & \leq  (n\pi)^{K_n/2} \cdot (1+R^2/n)^{-1/(16\mu_{n}^2)+K_n/2+1} / \Gamma(K_n/2), \end{align*}
for any $R\geq 0$ where, here, $\Gamma$ stands for the gamma function.
\end{lem}
\noindent Lemma \ref{L:FouTails} yields that $T^*$ has a continuous and bounded $K_n$-dimensional conditional Lebesgue density given $C$ whenever $\mu_n^{-2} > 8 K_n$. This density which we call $f_{n,\theta|C}^*$ turns out to be
\begin{equation} \label{eq:FouInv} f_{n,\theta|C}^*(x_1,\ldots,x_{K_n}) \, = \, (2\pi)^{-K_n} \idotsint \exp\Big(-i \sum_{k=1}^{K_n} t_k x_k\Big) \, \Psi_n^*(t_1,\ldots,t_{K_n}) \, dt_1\cdots dt_{K_n}\,, \end{equation}
for all $x=(x_1,\ldots,x_{K_n}) \in \mathbb{R}^{K_n}$ by Fourier inversion. The conditions of Lemma \ref{L:FouTails} are satisfied whenever $K_n\geq 2$, (\ref{eq:cond.pilot1}) and (\ref{eq:cond4'}) hold true.

\subsubsection{High-dimensional Edgeworth Expansion}

We need more precise investigation of the series representation of $\Psi_n^*$ in (\ref{eq:psin*}). In particular the remainder term is approximated by some $K_n$-variate polynomial with some degree not larger than $Q$; and an upper bound is needed for the coefficients and the new remainder term. This is provided in the following lemma.
\begin{lem} \label{L:Edge}
The remainder term in (\ref{eq:psin*}) admits the decomposition
\begin{align*}
\exp\Big\{\frac12 \sum_{\ell=3}^\infty \frac1\ell (2i)^\ell \cdot \mbox{tr}\big[\big(\tilde{D}_\theta^{[t]}\big)^\ell\big]\Big\} & \, = \, P_{Q,n}(t) \, + \, R_{Q,n}(t)\,, \end{align*}
where $P_{Q,n}$ denotes the $K_n$-variate polynomial
$$ P_{Q,n}(t) \, := \,  1 \, + \, \sum_{q=3}^Q \sum_{\|{\bf m}\|_1=q} \nu_{\bf m} \cdot t_1^{m_1} \cdots t_{K_n}^{m_{K_n}}\,, $$
with the degree $\leq Q$; where
$$ |\nu_{\bf m}| \, \leq \, {q \choose m_1,\ldots,m_{K_n}} \cdot 2^q\cdot \mu_n^{q/3}\cdot q^{1/4}\,, $$
when $\|{\bf m}\|_1 = \sum_{k=1}^{K_n} m_k = q$; and
\begin{align} \nonumber
|R_{Q,n}(t)| & \, \leq \, (Q+1)^{1/4} \cdot \big(2 K_n^{1/2} \mu_n^{1/3} \cdot \|t\| \big)^{Q+1} / \big\{1 - 2 K_n^{1/2} \mu_n^{1/3} (Q+1)^{1/(4Q+4)} \cdot \|t\|\big\}\,,
\end{align}
whenever $Q\geq 2$ and $\|t\| < (Q+1)^{1/(4Q+4)} \cdot\mu_n^{-1/3} K_n^{-1/2} /2$.
\end{lem}
\noindent By this decomposition we can prove asymptotic equivalence under polynomial growth of $K_n$ in $n$. Note that the constraints from Lemma \ref{L:Edge} also imply the weaker condition $\|t\| \leq 1/(2\mu_n)$ that ensures convergence of the expansion.

\subsubsection{Piecing Together the Ingredients for Gaussianization }

In the following, all probabilities and expectations shall be viewed as conditional given $C$. Note that the total variation distance between ${\cal L}_\theta(T^*|C)$ and ${\cal N}(0,I_{K_n})$ may be treated as follows.
\begin{align*} \sup_\theta & \, \mbox{TV}\big({\cal L}_\theta(T^*|C) \, , \, {\cal N}(0,I_{K_n})\big) \\ & \, = \, 1 \, - \, \inf_\theta \, \int \min\big\{f^*_{n,\theta|C}(x) , (2\pi)^{-K_n/2} \exp\big( - \|x\|^2/2\big)\big\} dx \\
& \, = \, 1 - \inf_\theta \, \mathbb{E}\, \min\big\{1 , f^*_{n,\theta|C}(X) \cdot (2\pi)^{K_n/2} \exp\big(\|X\|^2/2\big)\big\} \\
& \, \leq \, 1 - \int_0^1 \inf_\theta \, \mathbb{P}\big[f^*_{n,\theta|C}(X) > s \cdot (2\pi)^{-K_n/2} \exp\big(-\|X\|^2/2\big)\big]\, ds\,,
\end{align*}
where $x=(x_1,\ldots,x_{K_n})$ and $X$ denotes a ${\cal N}(0,I_{K_n})$-distributed random vector. By dominated convergence, the TV distance converges to zero uniformly with respect to $\theta$ if, for any fixed $s\in [0,1)$, it holds that
\begin{equation} \label{eq:TVmin} \liminf_{n\to\infty} \, \inf_\theta \, \mathbb{P} \big[f^*_{n,\theta|C}(X) > s \cdot (2\pi)^{-K_n/2} \exp\big(-\|X\|^2/2\big)\big] \, = \, 1\,. \end{equation}
We decompose the Fourier integral in (\ref{eq:FouInv}) as follows,
\begin{align*}
f^*_{n,\theta|C}(x) & \, \geq \, (2\pi)^{-K_n/2} \, \exp\big(-\|x\|^2/2\big) \, + \, \tau^{[1]}_{n,\theta}(x)  \, - \, \tau^{[2]}_{n,\theta}  \, - \, \tau^{[3]}_{n,\theta}  \, - \, \tau^{[4]}_{n,\theta}\,,
\end{align*}
where
\begin{align*}
\tau^{[1]}_{n,\theta}(x) & \, := \, \mbox{Re} \, (2\pi)^{-K_n} \, \int \exp\big( - i \langle x,t\rangle\big) \, \exp\big(-\|t\|^2/2\big) \cdot \big\{P_{Q,n}(t) - 1\big\} \, dt\,, \\
\tau^{[2]}_{n,\theta} & \, := \, (2\pi)^{-K_n} \, \int_{\|t\|\geq R_n} \exp\big(-\|t\|^2/2\big) \cdot \big|P_{Q,n}(t)\big| \, dt\,, \\
\tau^{[3]}_{n,\theta} & \, := \, (2\pi)^{-K_n} \, \int_{\|t\|< R_n} \exp\big(-\|t\|^2/2\big) \cdot \big|R_{Q,n}(t)\big| \, dt\,, \\
\tau^{[4]}_{n,\theta} & \, := \, (2\pi)^{-K_n} \, \int_{\|t\|\geq R_n} |\Psi_n^*(t)| \, dt\,,
\end{align*}
while $t=(t_1,\ldots,t_{K_n})$ and $R=R_n$ and $Q=Q_n$ remain to be selected. We study the asymptotic behaviour of these terms in the following lemma.
\begin{lem} \label{L:piecing}
    Under conditions (\ref{eq:cond.pilot1}), (\ref{eq:cond4'}) and
\begin{equation} \label{eq:Qn}
\lim_{n\to\infty} \, Q_n \, K_n^{4/3} n^{-1/3}\, = \, 0\,,
\end{equation}
we have
$$
\sup_\theta \mathbb{P}\big[\tau^{[1]}_{n,\theta}(X) \, \leq \, (s-1) \, (2\pi)^{-K_n/2} \exp\big(-\|X\|^2/2\big) / 2\big] \to 0
$$
as $n \to \infty$. Moreover, under the constraint that $R_n^2 \, \geq \, Q_n+K_n+1\,,$ if
\begin{align} \label{eq:condCLT4}
& R_n^2/n \, \longrightarrow \, 0\quad \text{and} \quad
R_n^2 K_n^{-2}\, / \, \log n \, \longrightarrow \, +\infty\,, \\
\label{eq:condCLT7}
& \frac12 R_n^2 K_n^{-1} \, - \, \log R_n \, \longrightarrow \, +\infty\,, \\
\label{eq:condCLT6}
 & K_n^4 \, R_n^6 / n \, \longrightarrow \, 0\,, \\
  \label{eq:condCLT8}
& Q_n/K_n \, \longrightarrow \, + \infty
\end{align}
as $n \to \infty$, we have
$$
\sup_\theta \mathbb{P}\big[ \tau^{[2]}_{n,\theta} \, + \, \tau^{[3]}_{n,\theta} \, + \, \tau^{[4]}_{n,\theta} \, \geq \, (1-s) \, (2\pi)^{-K_n/2} \exp\big(-\|X\|^2/2\big) / 2\big] \to 0.
$$
Thus, (\ref{eq:TVmin}) holds.
\end{lem}

Now we are ready to provide asymptotic equivalence of the experiments ${\cal F}_n$ and ${\cal G}_n$ and, hence, of ${\cal A}_n$ and ${\cal G}_n$ when we take into account the results of the previous sections. First, impose that
$$ \lim_{n\to\infty} \, K_n^{10} (\log n) / n \, = \, 0\,. $$
Then (\ref{eq:cond.pilot1}) and (\ref{eq:cond4'}) can be fulfilled by having $(\gamma_n)_n$ tend to $\infty$ just slightly faster than $(K_n)_n$. Next, we select $(Q_n)_n$ such that it tends to $\infty$ slightly faster than $(K_n)_n$ so that (\ref{eq:Qn}) and (\ref{eq:condCLT8}) hold true. Finally, choose $(R_n)_n$ such that it diverges to $\infty$ slightly faster than $\{K_n (\log n)^{1/2}\}_n$ under the constraint $R_n^2 \geq Q_n+K_n+1$ and satisfies (\ref{eq:condCLT4}), (\ref{eq:condCLT7}) and (\ref{eq:condCLT6}).

\begin{theo} \label{T.1}
Assume that $K_n\geq 2$ and
\begin{equation} \label{eq:Theo1}
\lim_{n\to\infty} \, K_n^{10} \, (\log n) \, / \, n \, = \, 0\,.
\end{equation}
The orthonormal system $(M_k)_k$ is chosen such that (\ref{eq:condMk}) is satisfied. Grant (\ref{eq:cond1}) and that all eigenvalues of all $\theta\in\Theta$ are bounded from below by $\rho>0$ and by $1/\rho$ from above. Then the experiments ${\cal A}_n$ and ${\cal G}_n$ are asymptotically equivalent.
\end{theo}
Let us note that the key element in this theorem is the convergence to 0 of the total variation distance between the distributions of the vector $T$ in \eqref{Tvec} and of the Gaussian vector $\mathcal{N}(d_\theta, \Gamma_\theta)$, conditionally on $C$. The vector $T$ is a multivariate functional of the standard Gaussian $n-$dimensional vector $\varepsilon$.

Many Gaussian approximations exist in the literature that can be obtained by Stein's formula and Malliavin calculus, see \cite{NourdinPeccatiBook}, and more generally by the method of moments and cumulants for Gaussian Wiener chaos in \cite{PeccatiTaqquBook}. The results show convergence in distribution, or use Wasserstein distances between the distributions. Bounds on the Wasserstein distance with explicit dependence on the dimension have been recently improved by \cite{Bonis2024}, see also references therein, in the CLT for a sum of independent vectors. 
In another line of work, it has been proven that convergence in Wasserstein distance can be turned into convergence in total variation distance under relatively mild assumptions, see \cite{BallyCaraPoly2020}, \cite{BallyCara2014} and references therein. More recently, it has been shown that for vectors with smooth probability density functions, the total variation distance can be bounded from above by some function of the Wasserstein $W_1$ distance, see \cite{Chae}. However, in these results there is no explicit control on how the constants depend on the dimension $K_n$ of the vector and, on closer inspection, the dimension appears in the exponent, whereas our result allows a polynomial growth of $K_n$, see assumption \eqref{eq:Theo1}.

\subsubsection{Transformation of the Covariance Matrices} \label{Trans}

Under the conditions of Lemma \ref{L:Gamma_theta}, which are satisfied under those of Theorem \ref{T.1}, we use (A.4) in \cite{Reiss2011} again to show that the experiment ${\cal G}_n$ is asymptotically equivalent to the experiment ${\cal H}_n$, in which $C$ and -- conditionally on $C$ -- a ${\cal N}(\Gamma \alpha_\theta/2,\Gamma)$-distributed random vector are observed, where $\alpha_\theta = \{\langle \theta,M_k\rangle_F\}_{k}^\top  =  \{\langle C_\theta,M_k\rangle_F\}_{k}^\top$. Therein we have applied that, for all $k=1,\ldots,K_n$,
\begin{align*}
d_{\theta,k} & \, = \, \mbox{tr}\big(C_\theta C^{-1} M_k C^{-1}\big) \, = \, \sum_{k'=1}^{K_n} \langle C_\theta,M_{k'}\rangle_F \cdot \mbox{tr}\big(M_{k'}C^{-1} M_k C^{-1}\big) 
\, = \, \frac12 \cdot \{\Gamma \alpha_\theta\}_k\,.
\end{align*}

By the action of the matrix $2 \Gamma^{-1}$ onto the second component observed under ${\cal H}_n$, this experiment is transformed to data from the equivalent experiment ${\cal I}_n$ in which $C$ and -- conditionally on $C$ -- a ${\cal N}(\alpha_\theta,4\Gamma^{-1})$-distributed random vector are observed, where the reverse transformation is provided by the action of $\Gamma/2$ to the second component of the observation from ${\cal I}_n$. Thus, we have arrived at a partial Gaussian experiment with a known covariance matrix (under the localization regime).

\subsubsection{Matrix Experiment with GOE Noise} \label{GOE}

We introduce the experiment ${\cal J}_n$ in which one observes $C$ from the experiment ${\cal I}_n$ and -- conditionally on $C$ -- the random matrix
$$ \hat{C} \, := \, C^{-1/2} C_\theta C^{-1/2} \, + \, G\,, $$
where $G$ denotes the $n$-dimensional Gaussian Orthogonal Ensemble (GOE), i.e. the random symmetric $n\times n$-matrix $G$ where all entries $G_{j,j'}$, for $j\geq j'$, are independent with $G_{j,j'} \sim {\cal N}(0,1)$, $j>j'$, and $G_{j,j}\sim {\cal N}(0,2)$, $j=1,\ldots,n$. Given $C$, the random matrix $\hat{C}$ has the conditional density
\begin{align*} f_{\hat{C},\theta}(X) & \, = \, 2^{-1/2} (2\pi)^{-n(n+1)/4} \cdot \exp\big\{ - \big\|X-C^{-1/2}C_\theta C^{-1/2}\big\|_F^2 / 4\big\} \\
& \, = \, 2^{-1/2} (2\pi)^{-n(n+1)/4} \cdot \exp\big\{ - \|X\|_F^2 / 4\big\} \cdot \exp\big\{ - \big\|C^{-1/2} C_\theta C^{-1/2}\big\|_F^2 / 4\big\} \\
 & \hspace{0.2cm} \cdot \exp\Big\{ - \sum_{k=1}^{K_n} \alpha_{\theta,k} \cdot \mbox{tr}\big(  C^{-1/2} M_k C^{-1/2} X\big) / 2\Big\}\,,\end{align*}
for all symmetric $n\times n$-matrices $X$ with respect to the $n(n+1)/2$-dimensional Lebesgue measure on the vectorized upper triangular matrix of $X$ (including the main diagonal). By the Fisher-Neyman factorization,  $C$ and $\big\{\mbox{tr}(C^{-1/2} M_k C^{-1/2} \hat{C})\big\}_{k=1,\ldots,K_n}$ form a sufficient statistic in the experiment ${\cal J}_n$. As
$$ \big\{\mbox{tr}(C^{-1/2} M_k C^{-1/2} \hat{C})\big\}_k \, = \, \Gamma \alpha_\theta/2 \, + \, \big\{\mbox{tr}(C^{-1/2} M_k C^{-1/2} G)\big\}_k $$
for all $k=1,\ldots,K_n$, the conditional distribution of $\big\{ \mbox{tr}(C^{-1/2} M_k C^{-1/2} \hat{C}\big\}_k$ given $C$ equals ${\cal N}(\Gamma \alpha_\theta/2,\Gamma)$, so that model ${\cal J}_n$ has been shown to be equivalent to the experiment ${\cal H}_n$. By adding $I_n$ to $\hat{C}$ and subtracting it for the reverse transformation, the experiment ${\cal K}_n$ which describes the observation of $C$ and $C^{-1/2} \Delta_\theta C^{-1/2} + G$ is equivalent to ${\cal J}_n$. Therein note that $-G$ and $G$ are identically distributed.

\subsection{Re-Globalization}

We introduce the statistical experiment ${\cal L}_n$ in which one observes a ${\cal N}(\alpha_\theta,4\tilde{\Gamma}_\theta^{-1})$-distributed random variable with the matrix $\tilde{\Gamma}_\theta$ from (\ref{eq:tildeGamma}) and the parameter set $\Theta$ as in the experiment ${\cal A}_n$. As a pilot estimator $\hat{\alpha}$ for $\alpha_\theta$ in ${\cal L}_n$ we propose the observation itself, which is unbiased and satisfies
\begin{align*}
\sup_\theta \, & \mathbb{E}_\theta \|\hat{\alpha} - \alpha_\theta\|^2 \, = \, 4\, \sup_\theta \, \mbox{tr}\big(\tilde{\Gamma}_\theta^{-1}\big) \, \leq \, 4 K_n \, \rho^{-2}
\end{align*}
for $n$ sufficiently large (uniformly with respect to  $\theta$), since, by Lemma \ref{L:Gamma_theta}, all eigenvalues of the symmetric matrix $\tilde{\Gamma}_\theta^{-1}$ are bounded from above by $\rho^{-2}$. Therefore we can use the localization scheme from Section \ref{GLS} again so that, by appropriate choice of the sequence $\beta_n$ and under the constraint (\ref{eq:cond.pilot1}), the experiment ${\cal L}_n$ is asymptotically equivalent to the experiment, in which -- in addition to and independently of the observation from ${\cal L}_n$ -- the matrix $C$ as in the experiment ${\cal I}_n$ is accessible. Then, under the conditions of Theorem \ref{T.1}, we combine Lemma~\ref{L:Gamma_theta} and
the following auxiliary result.

\begin{lem} \label{L:sp.1}
Let $A$ and $B$ be symmetric, positive semi-definite matrices of the same size, where $A$ is assumed to be even positive definite. If $\|A^{-1/2}(B-A)A^{-1/2}\|_{sp} < 1$, then $B$ is positive definite and
$$ \big\|A^{1/2}(A^{-1}-B^{-1})A^{1/2}\big\|_F \, \leq \, \|A^{-1/2}(B-A)A^{-1/2}\|_F \, \big/ \, \big(1 - \|A^{-1/2}(B-A)A^{-1/2}\|_{sp}\big)\,.  $$
\end{lem}
This proves asymptotic equivalence of the experiments ${\cal L}_n$ and ${\cal I}_n$ so that we arrive at the following theorem, which provides asymptotic equivalence of two global experiments, i.e. experiments without pre-localized parameters.
\begin{theo} \label{T.2}
Under the conditions of Theorem \ref{T.1}, the experiments ${\cal A}_n$, ${\cal H}_n$, ${\cal I}_n$, ${\cal J}_n$, ${\cal K}_n$ and ${\cal L}_n$ are asymptotically equivalent.
\end{theo}

\section{Asymptotic Equivalence of Locally Stationary Times Series and Bivariate White Noise} \label{sec:AE_LSP_BGW}

Now we consider data which are drawn from a centered Gaussian locally stationary time series; and apply the general scheme from Section \ref{GGS} to establish asymptotic equivalence of this model and of the bivariate Gaussian white noise model.

\subsection{Locally Stationary Time Series} \label{LocStat}
In the sense of \cite{Dahlhaus1996,Dahlhaus1997,10.1214/aos/1015957480}, a triangular array $(X_{t,n})_{1\leq t\leq n}$ of centered Gaussian processes (in discrete time) is locally stationary with transfer function $A^0=(A^0_{t,n})_{1\leq t\leq n}$,
\begin{itemize}
    \item[(i)]
if it possesses a time-varying Cram{\'e}r type representation
$$
X_{t,n}=\int_{-\pi}^{\pi}\exp(ix t)A^0_{t,n}(x)\,dZ(x),\ \ t=1,\dots, n,\, n\in\mathbb{N},
$$
where $Z(-x)=\overline{Z(x)}$ and $Z$ is restricted to $[0,\pi]$ a complex Brownian motion, and
\item[(ii)] if there exists a function  $A:[0,1]\times\mathbb{R}\rightarrow\mathbb{C}$ with $A(u,-x)=\overline{A(u,x)}$ that is continuous in its first argument,  $2\pi$-periodic in its second argument and satisfies
$$
\sup_n \Big(n\cdot \sup_{t,x}\big\arrowvert A^0_{t,n}(x)-A(t/n,x)\big\arrowvert\Big)<\infty.
$$
\end{itemize}
The function $f:[0,1]\times |-\pi,\pi]\rightarrow\mathbb{R}$ with $f(u,x)=A(u,x)\overline{A(u,x)}$ is called time-varying spectral density. By (ii), it is symmetric in its second argument.

\begin{rem}
By mirroring $A$ and  $A^0$ at the vertical axis through $(1,0)$, we may extend them to $A^0_{1,n},\dots, A^0_{n,n},A^0_{t+1,n},\dots,  A^0_{2n,n}$ and
$A:[0,2]\times[-\pi,\pi] \rightarrow \mathbb{C}$,  such that the spectral density $f:[0,2]\times [-\pi,\pi]\rightarrow\mathbb{R}$ is periodic also in its first argument. As the upper left quarter of the enlarged  $(2n\times 2n)$-covariance matrix corresponds to the original ($n\times n)$-matrix and the extended spectral density shares the same upper and lower bounds, we assume subsequently that $f(0,\cdot)=f(1,\cdot )$.
\end{rem}
Define the functions of the Fourier basis $\tilde{\varphi}_{j,j'}$, for $j,j' \in \mathbb{Z}$, by
$$
\tilde{\varphi}_{j,j'}(t,x) \, := \, \exp\big(2\pi i j t + i j' x\big), \quad t\in [0,1], \, x \in [-\pi,\pi]\,.
$$
They form an orthonormal basis of the square-integrable complex-valued functions defined on $[0,1]\times [-\pi, \pi]$. Such a function $f$ belongs to a Sobolev ellipsoid $S^2(s,\tilde L)$ with global smoothness $s>0$ and radius $\tilde L>0$ if
$$
\sum_{j,j' } |\langle f, \tilde \varphi_{j,j'} \rangle_{L_2}|^2 [1+ (j^2 + (j')^2)^s] \leq \tilde L.
$$
In order to take into account that spectral density functions $f$ are real-valued and symmetric with respect to the frequency argument $x$ on $[-\pi,\pi]$, we introduce the following orthonormal system of functions
\begin{equation} \label{eq:varphi+}
\varphi_{j,j'}^+(t,x) \, := \, \begin{cases} \sqrt{2/\pi} \, \cos(j'x) \cos(2\pi j t)\,, & \mbox{ if }j,j' \neq 0\,, \\
\sqrt{1/\pi} \, \cos(j'x)\,, & \mbox{ if }j=0, j'\neq 0\,, \\
\sqrt{1/\pi} \, \cos(2\pi j t)\,, & \mbox{ if }j\neq 0, j'= 0\,, \\
\sqrt{1/(2\pi)} \,, & \mbox{ otherwise,} \end{cases} \end{equation}
for $j,j'\geq 0$; and, for $j\geq 1, j'\geq 0$,
\begin{equation} \label{eq:varphi-} \varphi_{j,j'}^-(t,x) \, := \, \begin{cases} \sqrt{2/\pi} \, \cos(j'x) \sin(2\pi j t)\,, & \mbox{ if }j' \neq 0\,, \\
\sqrt{1/\pi} \, \sin(2\pi j t)\,, & \mbox{ if }j'= 0\,,
\end{cases} \end{equation}
which form an orthonormal system in the Hilbert space $L_2([0,1]\times[-\pi,\pi])$ of all measurable squared-integrable functions on the domain $[0,1]\times[-\pi,\pi]$. Moreover it is an orthonormal basis of the subspace of those elements which are real-valued and symmetric in the second component. The corresponding inner product is denoted by $\langle\cdot,\cdot\rangle_{L_2}$. Besides being real-valued and symmetric in its second argument,
we assume that the function $f$  belongs to the Sobolev ellipsoid with global smoothness $s>0$ and radius $L>0$ which is characterized by the property
$$
\sum_{j,j'\geq 0} \big(|\tilde{f}_{j,j'}^+|^2 + |\tilde{f}_{j,j'}^-|^2\big) \, \big(j^2 + (j')^2\big)^s  \leq L,
$$
where $\tilde f_{j,j'}^\pm = \langle f, \varphi_{j,j'}^\pm \rangle_{L_2}$ denote the Fourier coefficients of the function $f$;  $\tilde{f}_{0,j'}^- :=0$ for convenience. Moreover, we assume that there exists some fixed constant $\rho^*$ in (0,1) such that $0<\rho^* \leq f(t,x) \leq 1/\rho^*$ for all $(t,x) \in [0,1] \times (-\pi,\pi)$. Let us denote the class of such spectral density functions by $W(s,L,\rho^*)$.


If the triangular array of Gaussian random variables $(X_{t,n})_{1\leq t\leq n}$ is locally stationary, it is immediate from the definition that
\begin{equation}\label{eq: vartheta}
\vartheta_{s,t}^n:=\mathbb{E}\big(X_{s,n}X_{t,n}\big)=\int_{-\pi}^{\pi}\exp(ix(s-t))A^0_{s,n}(x)\overline{A^0_{t,n}(x)}dx.
\end{equation}
As the targeted bivariate Gaussian white noise model is fully described by the time-varying spectral density $f$, we shall express the covariance matrices in the experiment $\mathcal{A}_n$ solely in terms thereof.  For this purpose, the covariance matrix $\theta = \{\theta(f)\}_{j,j'=1,\ldots,n}$ in the experiment ${\cal A}_n$ is assumed to have the shape
\begin{equation} \label{eq:density}
\theta_{j+1,j'+1}(f) \, = \, \int_{-\pi}^\pi \exp\big(i(j-j')x\big) \, f(\min\{j,j'\}/n,x) \, dx\,,
\end{equation}
for all $j,j'=0,\ldots,n-1$.  The next lemma justifies this approximation of the covariance matrix $\vartheta^n=(\vartheta^n_{s,t})_{1\leq s,t\leq n}$ in \eqref{eq: vartheta}. It reveals  that the centered $n$-dimensional normal vectors with covariances $\vartheta^n$ and $\theta(f)$ approximate each other in Hellinger distance as $n\rightarrow\infty$ under some mild smoothness constraint. We set $A_{t,n}^0(u,x):=A_{t,n}^0(x)$ for $u\in[0,1]$.
\begin{lem}\label{lem: 3.1} Let $s>5/2$ and $\tilde L,  L>0$. If $A, A_{t,n}^0\in S^2(s,\tilde L)$ for  $t=1,\dots, n$ and $n\in\mathbb{N}$, then we have
    $
\Arrowvert \vartheta^n-\theta(f) \Arrowvert_F=o(1)$. Moreover, $\Arrowvert \theta(f)\Arrowvert_{sp}\leq 2\pi/\rho^* + o(1)$, $\Arrowvert \theta(f)^{-1}\Arrowvert_{sp}\leq (2\pi/\rho^*+o(1))^{-1} $ for $f\in W(s,L,\rho^*)$.
\end{lem}

 Note that the functions $\varphi^+_{k,k'}$ and $\varphi^-_{k,k'}$ are obtained via the following linear transformations from the complex-valued Fourier basis $\tilde \varphi_{k,k'}$:
\begin{align} \label{eq:tildevarphi+}
\varphi_{k,k'}^+ & \, = \, \big(\tilde{\varphi}_{k,k'} + \tilde{\varphi}_{-k,k'} + \tilde{\varphi}_{k,-k'} + \tilde{\varphi}_{-k,-k'}\big) \cdot \frac14 \cdot  \begin{cases} \sqrt{2/\pi}\,, & \mbox{ if }k,k' \neq 0\,, \\
\sqrt{1/\pi} \,, & \mbox{ if }k=0, k'\neq 0\,, \\
\sqrt{1/\pi}\,, & \mbox{ if }k\neq 0, k'= 0\,, \\
\sqrt{1/(2\pi)} \,, & \mbox{ otherwise,} \end{cases} \end{align}
and
\begin{align} \label{eq:tildevarphi-}
\varphi_{k,k'}^- & \, = \, \big(\tilde{\varphi}_{k,k'} - \tilde{\varphi}_{-k,k'} + \tilde{\varphi}_{k,-k'} - \tilde{\varphi}_{-k,-k'}\big) \cdot \frac1{4i} \cdot  \begin{cases}  \sqrt{2/\pi} \,, & \mbox{ if }k' \neq 0\,, \\
\sqrt{1/\pi} \,, & \mbox{ if }k'= 0 \end{cases}\,.
\end{align}
It is therefore obvious that the projection coefficients on these basis functions check the same linear relations and thus any spectral density function $f$ in $W(s,L,\rho^*)$ can be shown to belong to a Sobolev ellipsoid with the same smoothness $s$ and some radius $\tilde L$ depending on $L$. Let fix some enumeration of the functions $\varphi_{j,j'}^+$, $j=0,\ldots,\kappa_1$, $j'=0,\ldots,\kappa_2$ and $\varphi_{j,j'}^-$, $j=1,\ldots,\kappa_1$, $j'=0,\ldots,\kappa_2$, from (\ref{eq:varphi+}) and (\ref{eq:varphi-}), with $\kappa_1,\kappa_2$ some positive integers. Denote by $\{\varphi_k : k=1,\ldots,K_n\}$ the enumeration of these basis functions. This enumeration is continued such that the $\varphi_k$, $k>K_n$, denote the $\varphi^\pm_{j,j'}$ with $j>\kappa_1$ or $j'>\kappa_2$.

\subsection{Useful Matrices for Locally Stationary Processes} \label{CircMat}

In the analysis of stationary processes, \cite{GolNussZhou2010} show the importance of circulant matrices in mapping Toeplitz covariance matrices into spectral density functions projected on Fourier series. In the case of locally stationary time series the covariance matrix is not Toeplitz, but the diagonals are associated to smooth functions which are well approximated by their projection on Fourier series. In Appendix A  in \cite{AELS_supplementary}, we construct  the matrices $\{\check{M}_{k}\}_{k=1,\ldots,K_n}$ which are mapped into the orthonormal system of functions
$\{\varphi_{k}\}_{k=1,\ldots,K_n}$ by a unique linear map $\Psi$.  When the codomain of $\Psi$ is restricted to its range it forms an isometric isomorphism between two Hilbert spaces. Moreover, $\Psi$ is shown to be an approximate homomorphism with respect to standard matrix multiplication   and pointwise multiplication of functions, respectively (Lemma A.1). In order to recover the Fourier coefficients $\langle f,\varphi_k \rangle_{L_2}$ of the spectral density $f$  in the projection coefficients of the covariance matrix $\theta(f)$ on some orthonormal system of matrices, we then build in Section \ref{S:BS} the matrices $\{M_{k}\}_{k=1,\ldots,K_n}$ and $\langle \theta(f), M_k\rangle_F$. The matrices $\{\check{M}_{k}\}_{k}$ and $\{M_{k}\}_{k}$ are shown to be close in Frobenius norm.
In Section \ref{subsec: 3.2.3}, sufficient conditions are stated for deriving asymptotic equivalence of the locally stationary time series and the  Gaussian experiment $\mathcal{L}_n$.

\subsubsection{Construction of the Basis $\{M_k\}_k$}\label{S:BS}

{Let $f^{[K_n]}$ denote the orthogonal projection of $f$ onto the linear hull of the $\varphi_{j,j'}^+$, $j=0,\ldots,\kappa_1$, $j'=0,\ldots,\kappa_2$ and the $\varphi_{j,j'}^-$, $j=1,\ldots,\kappa_1$, $j'=0,\ldots,\kappa_2$, from (\ref{eq:varphi+}) and (\ref{eq:varphi-}), with $\kappa_1,\kappa_2$ as in the supplemental material. Without loss of generality, we assume we adopt the same enumeration of the basis functions $\{\varphi_k : k=1,\ldots,K_n\}$ and for the matrices $\check{M}_k$. This enumeration is continued such that the $\varphi_k$, $k>K_n$, denote the $\varphi^\pm_{j,j'}$ with $j>\kappa_1$ or $j'>\kappa_2$. The linearity of the map $f \mapsto \theta(f)$ yields that
\begin{align*}
\theta\big(f^{[K_n]}\big) & \, = \, \sum_{k=1}^{K_n} \langle f,\varphi_k\rangle_{L_2} \cdot \theta(\varphi_k)\,. \end{align*}
Since
$$ \int_{-\pi}^\pi \exp(i(j-j')x)\, \tilde{\varphi}_{k,k'}(\min\{j,j'\}/n,x) \, dx \, = \, 2\pi \cdot \lambda\big(k \cdot \min\{j,j'\}/n\big) \cdot {\bf 1}\{j'=j+k'\}\,, $$
we derive that $\theta\big(\tilde{\varphi}_{k,k'}\big) \, = \, 2\pi \, \Lambda^k M^{k'}$ with $\tilde{\varphi}_{k,k'}$, $\lambda$ and $\Lambda$ as in Appendix A; and
$$
M \, := \, \begin{pmatrix} 0 & 1 & 0 & \cdots & 0 \\
                              \vdots &  &  & \ddots & \vdots \\
                              0 & 0 & 0 & \cdots & 1 \\
                              0 & 0 & 0 & \cdots & 0 \end{pmatrix}\,,
$$
where, for $k'<0$, we put $M^{k'} := (M^\top)^{-k'}$. Moreover, define
$$
\lambda^*_{j,j'}:= \exp(2\pi i j/(n-|j'|))
$$ for all $|j'|\leq n-1$ and $j$ such that $|j| <  n-|j'|$ and the diagonal matrix $\Lambda_{j,j'}$ having main diagonal entries $\big(\lambda_{j,j'}^0, \ldots , \lambda_{j,j'}^{n-|j'|-1}, 0,\ldots,0\big)$. Let us define $M_{0,j,j'} = 2\pi \Lambda_{j,j'} M^{j'}$. We deduce that
\begin{align*}
\big\|\theta\big(\tilde{\varphi}_{k,k'}\big) - M_{0,k,k'}\big\|_F^2
&  \leq  \sum_{j=0}^{n-k'-1} 16\pi^4 k^2 j^2 \cdot \{1/n - 1/(n-k')\}^2  \leq  16\pi^4 \, \kappa_1^2 \kappa_2^2 / n \asymp K_n^2 / n,
\end{align*}
under the constraints of Lemma \ref{L:matrixf}(a), which are imposed in the following.

We construct the matrices $M^\pm_{k,k'}$ by replacing $\tilde{\varphi}_{\pm k,\pm k'}$ by $M_{0,\pm k,\pm k'}$ on the right side of the equations (\ref{eq:tildevarphi+}) and (\ref{eq:tildevarphi-}). Then,
\begin{align} \label{eq:mkk+-}
\big\|\theta\big(\varphi^\pm_{k,k'}\big) - M^\pm_{k,k'}\big\|_F^2 & \, \leq \,  32\pi^3 \, \kappa_1^2 \kappa_2^2 / n \, \asymp\, K_n^2/n\,.
\end{align}
By $M_k^*$ we denote that $M^\pm_{j,j'}$ which corresponds to $\varphi_k = \varphi^\pm_{j,j'}$ in the arranged enumeration.
Moreover,
\begin{align*} & \langle M_{0,j,j'} , M_{0,k,k'} \rangle_F \, = \, 4\pi^2 \cdot \mbox{tr}\big(\Lambda_{-k,k'} \Lambda_{j,j'} M^{j'-k'}\big) \\ & \, = \, 4\pi^2 \cdot {\bf 1}\{j'=k'\} \, \sum_{\ell=0}^{n-|k'|-1} \exp\Big(2\pi i \ell \frac{j-k}{n-|k'|}\Big)  \, = \, 4\pi^2 (n-|k'|) \cdot {\bf 1}\{j'=k'\} \cdot {\bf 1}\{j=k\}\,, \end{align*}
for all $|j|,|k|\leq \kappa_1$, $|j'|,|k'| \leq \kappa_2$ with $\kappa_1,\kappa_2 < n/4$. Considering (\ref{eq:tildevarphi+}) and (\ref{eq:tildevarphi-}) it follows that the $ M_{j,j'}^+$, $j=0,\ldots,\kappa_1$, $j'=0,\ldots,\kappa_2$, and the
$M_{j,j'}^-$, $j=1,\ldots,\kappa_1$, $j'=0,\ldots,\kappa_2$ form an orthogonal system and satisfy
$$
\|M_{j,j'}^\pm\|_F^2 = 2\pi(n-j').$$
We arrive at the orthonormal system
$$
\{ M_k, \,k=1,\ldots,K_n \}, \quad M_k := M_k^* / \|M_k^*\|_F, \text{ where }2\pi n \geq \|M_k^*\|_F^2 \geq 2\pi(n-\kappa_2).
$$

Finally we show the proximity between the matrices $M_k$ and $\check{M}_k$. Note that
$$ \big\|\theta\big(\tilde{\varphi}_{k,k'}\big) - 2\pi \check{M}_{0,k,k'}\big\|_F^2
 = 4 \pi^2 \big\|\Lambda^k \cdot (M^{k'} - \check{M}^{k'})\big\|_F^2\, = \, 4\pi^2 |k'|\,, $$
from what follows that
$$ \big\|\theta\big(\varphi^\pm_{k,k'}\big) - 2\pi \check{M}^\pm_{k,k'}\big\|_F^2 \, \leq \, 8\pi |k'|\,, $$
by (\ref{eq:varphi+}) and (\ref{eq:varphi-}). Combining this result with (\ref{eq:mkk+-}) yields that
$$ \sup_{j\leq \kappa_1,j'\leq \kappa_2} \big\|M_{j,j'}^\pm - 2\pi \check{M}_{j,j'}^\pm\big\|_F^2 \, = \, {\cal O}\big(K_n^{1/2} + K_n^2 n^{-1}\big)\,. $$
For the corresponding normalized matrices we obtain
\begin{align*} \bigg\|\frac{M_{j,j'}^\pm}{\|M_{j,j'}^\pm\|_F} - \frac{\check{M}_{j,j'}^\pm}{\|\check{M}_{j,j'}^\pm\|_F}\bigg\|_F &  \leq \frac1{\sqrt{2\pi n}}  \|2\pi \check{M}_{j,j'}^\pm - M_{j,j'}^\pm\|_F  + \frac{\|M_{j,j'}^\pm\|_F}{\sqrt{2\pi}}  \Big|\frac1{\sqrt{n-j'}} - \frac1{\sqrt{n}}\Big|, \end{align*}
so that
\begin{equation} \label{eq:MMcheck}
\max_{k=1,\ldots,K_n} \|M_k - \check{M}_k\|_F \, = \, {\cal O}\big(K_n^{1/4} n^{-1/2} + K_n n^{-1} + K_n^{1/2} n^{-1}\big)\,, \end{equation}
where, under (\ref{eq:Theo1}), the term $K_n^{1/4} n^{-1/2}$ dominates the others.}

\subsubsection{Sufficient conditions for asymptotic equivalence}\label{subsec: 3.2.3}

We show the following results.
\begin{lem} \label{L:matrixf}
(a)\; We choose $\kappa_1$ and $\kappa_2$ such that $\kappa_1 \asymp \kappa_2 \asymp \sqrt{K_n}$. Then, for $s>1$,
\begin{align*} & \sup_{f\in W(s,L,\rho^*)} \|f-f^{[K_n]}\|_\infty \, = \, {\cal O}\big(K_n^{(1-s)/2}\big) \\ \mbox{and }&
\sup_{f\in W(s,L,\rho^*)} \|f-f^{[K_n]}\|_{L_2} \, = \, {\cal O}\big(K_n^{-s/2}\big)\,.
\end{align*}
\noindent (b)\; Let $f,\tilde{f}\in  W(s,L,\rho^*)$ be real-valued and symmetric in the second argument. Then,
\begin{align*}
    \|\theta(f) - \theta(\tilde f)\|_F^2 \leq 12 \pi^2 \cdot n \cdot \|f - \tilde f\|_\infty^2   \, .
\end{align*}
There exist positive constants $C_0$ and $C_1$ depending only on $s$ and $L$ such that
$$
\sup_{f\in W(s,L,\rho^*)} \|f\|_\infty \leq C_0 \quad \text{if } s>1 \text{ and }
\sup_{f\in W(s,L,\rho^*)} \|\frac{\partial}{\partial t}f (t,x)\|_\infty \leq C_1 \quad \text{if } s>2.
$$
Moreover, for $s>2$, we have:
$$\|\theta(f)-\theta(\tilde{f})\|_F^2 \, \leq \, 6 \pi  \cdot n \cdot \left( \|f-\tilde{f}\|_{L_2}^2 + \frac{4 \pi C_1 \|f - \tilde f\|_\infty}{n}\right)\,.
$$
\end{lem}
 Let us apply Lemma~\ref{L:matrixf} to $\tilde f = f^{[K_n]}$. We get that
 $$
 \sup_{f\in W(s,L,\rho^*)} \|\theta(f) - \theta(f^{[K_n]})\|_F^2 = {\cal O} (n K_n^{-s} + K_n^{(1-s)/2} ).
 $$ Moreover, by the Cauchy-Schwarz inequality,
\begin{align*}
\sup_{f\in W(s,L,\rho^*)} \Big\|\theta\big(f^{[K_n]}\big) - \sum_{k=1}^{K_n} \langle f,\varphi_k\rangle_{L_2}\cdot M^*_k\Big\|_F^2
&  \leq \sup_{f\in W(s,L,\rho^*)} \|f\|^2_{L_2} \cdot \sum_{k=1}^{K_n} \|\theta(\varphi_k) - M_k^*\|_F^2 \\
& = {\cal O}(K_n^3 /n)\,,
\end{align*}
since $\|f\|_{L_2}^2 \leq 2\pi/ (\rho^*)^{2}$. Combining this inequality with Lemma \ref{L:matrixf} provides that
\begin{equation} \label{eq:Ctheta}
\sup_{f\in W(s,L,\rho^*)} \, \Big\| \theta(f) - \sum_{k=1}^{K_n} \langle f,\varphi_k \rangle_{L_2}\cdot  M_k^*\Big\|_F \, = \, {\cal O}\big( n^{1/2}  K_n^{-s/2} + K_n^{(1-s)/4}  + K_n^{3/2} n^{-1/2} \big)\,.	
\end{equation}

Note that
$$ C_{\theta(f)} \, = \,  \sum_{k=1}^{K_n} \alpha_{\theta(f),k}\cdot M_k \, = \, \sum_{k=1}^{K_n} \langle f,\varphi_k\rangle_{L_2}\cdot M_k^*\,, $$
where $C_{\theta(f)}$ is as in Section \ref{PreSmooth} and
\begin{equation} \label{eq:coeffalpha} \alpha_{\theta(f),k} = \langle \theta(f),M_k\rangle_F = \langle f,\varphi_k\rangle_{L_2} \cdot \|M_k^*\|_F. \end{equation}
 Since $\|M_k^*\|_{sp} \leq 2\sqrt{2\pi}$ for all $k=1,\ldots,K_n$ we verify condition (\ref{eq:condMk}). Furthermore,
\begin{align*}
\big\|\theta(f)^{-1/2} \big(\theta(f) - C_{\theta(f)}\big) \theta(f)^{-1/2}\big\|_F^2 & \, \leq \, \big\|\theta(f)^{-1/2}\big\|_{sp}^4 \cdot \big\|\theta(f) - C_{\theta(f)}\big\|_F^2 \\ & \, \leq \, \rho^{-2} \cdot \big\|\theta(f) - C_{\theta(f)}\big\|_F^2\,,
\end{align*}
so that, by (\ref{eq:Ctheta}), the condition (\ref{eq:cond1}) is satisfied whenever
\begin{equation} \label{eq:condspec1}
\lim_{n\to\infty} \,  n^{1/2}  K_n^{-s/2} + K_n^{(1-s)/4}  + K_n^{3/2} n^{-1/2}  \, = \, 0\,.
\end{equation}
Hence, Theorems \ref{T.1} and \ref{T.2} are applicable whenever $s>10$. This condition enables us to select $K_n$ such that (\ref{eq:Theo1}) and (\ref{eq:condspec1}) are fulfilled and to state the following result.

\begin{theo} \label{AE_LSP}
Let $\mathcal{A}_n$ be the experiment of a Gaussian locally stationary time series
$$
\left( \mathbb{R}^n, \mathcal{B}(\mathbb{R}^n), \{ \mathcal{N}(0, \theta(f)): f \in W(s,L,\rho^*)\}\right)
$$
with $s>10$ and let the matrices $\{M_k; k = 1,\ldots,K_n\}$ be defined in Section~\ref{S:BS}.  If $K_n$ is such that $K_n^{10} \log n /n \to 0$ as $n \to \infty$, then the experiments $\mathcal{A}_n$, $\mathcal{K}_n$ and $\mathcal{L}_n$ are equivalent.
\end{theo}

\subsection{Bivariate Gaussian White Noise} \label{LSDM}

Let us consider the experiment ${\cal A}^*_n$ where we observe the bivariate Gaussian process $X(t,x)$, $t\in[0,1]$ and $x \in (-\pi,\pi)$ defined by
$$
dX(t,x) = \log f(t,x) dt \, dx + 2 \sqrt{\frac{\pi}{n}} \cdot dW(t,x),
$$
where $W$ is a Wiener process on $\mathbb{R}^2$ and $f$ belongs to the class $W(s,L,\rho^*)$.
In the rest of the section we follow the general methodology of  Section \ref{GGS} for Gaussianization of this model. First, we propose a pilot estimator, then we give sufficient conditions such that the general scheme applies.

\subsubsection{Pilot Estimator}

The first step is to estimate the drift $\log f $ in this model by projecting on the orthogonal functions  $\varphi_j$ for $j=1,\ldots, J_n$ from Section \ref{LocStat}:
$$
\widetilde{\log } f (t,x) \, = \, \frac1{\sqrt{2\pi n}} \, \sum_{j=1}^{J_n}  \tilde \alpha_j \cdot \varphi_j(t,x) ,
$$
where $\tilde \alpha_j := \sqrt{2\pi n} \cdot \langle \varphi_j , dX \rangle_{L_2} $.
Next, we build an estimator of $f$:
$$
\hat f(t,x) \, = \, \frac1{\sqrt{2\pi n}} \, \sum_{j=1}^{K_n} \hat \alpha_j \cdot \varphi_j(t,x)\,,
$$
where $\alpha_j := \sqrt{2 \pi n} \cdot \langle \exp( \widetilde{\log} f ), \varphi_j \rangle_{L_2}$ and show sufficient conditions such that $\alpha_j $ can be written under the form $ \alpha_{\theta(f),j} + \tilde \eta_j$ with $\|\tilde \eta\| \leq \gamma_n$, $\alpha_{\theta(f),j}$ from \eqref{eq:coeffalpha} which can also be written as $\alpha_{\theta(f),j} = a_j \|M_j^*\|_F / \sqrt{2 \pi n}$ for $a_j = \sqrt{2 \pi n} \cdot \langle f, \varphi_j \rangle_{L_2}$.

\begin{lem} \label{lem:pilot} In the model $\mathcal{A}_n^*$ with $s\geq 3$, the estimator $\hat f$, or the $K_n-$dimensional vector $ \alpha$, is such that
$$
\sup_{f \in W(s,L,\rho^*)} \mathbb{E} \left[ \| \hat \alpha - \alpha_{\theta(f)}\|^2 \right] \leq \mathcal{O}(1) \cdot K_n \, .
$$
Thus, \eqref{eq:cond2} and \eqref{eq:cond3} are satisfied, provided that $K_n/n^2 = \mathcal{O}(1)$ and that $K_n/\gamma_n \to 0$.
\end{lem}

\subsubsection{Localized Gaussian Bivariate White Noise}

We impose that
$$
f_n(t,x) \, = \, \frac1{\sqrt{2\pi n}} \, \sum_{j=1}^{K_n} \, \alpha_{\theta(f),j} \cdot \varphi_{j}(t,x)\,,
$$
with the functions $\varphi_j$ from Section \ref{LocStat}. Moreover, based on the previous section,  we impose that
$$
\hat{f}(t,x) \, = \, \frac1{\sqrt{2\pi n}} \, \sum_{j=1}^{K_n} \hat \alpha_j \cdot \varphi_j(t,x)\,,
$$
where the $\hat \alpha_j := \alpha_{\theta(f),j} + \tilde{\eta}_j$ are gained from some independent additional information on $\alpha_\theta$ (see Section \ref{GLS}). The localized ${\cal A}_n^*$-model $(\mathcal{A}_n^*)'$ is further equivalent to observing $\hat f$ and
$$
dX_n(t,x) = \log f_n(t,x) \, dt \, dx + 2\sqrt{\frac \pi n}\,  dW(t,x),
$$
where $W$ denotes a standard Wiener process on $\mathbb{R}^2$ (two-dimensional Brownian sheet).  Indeed, it is sufficient to prove that
\begin{align} \label{equiv1}
    \frac n{4 \pi} \|\log f - \log (f_n) \|_{L_2}^2 \to 0, \quad n \to \infty.
\end{align}
Let us first prove that $f_n$ and $\hat f$ belong to $[\rho^*/2, 2/ \rho^*]$.
%
By $\|\varphi_k\|_\infty \leq \sqrt{2/\pi}$, (\ref{eq:coeffalpha}), the Cauchy-Schwarz inequality, we derive that
\begin{align*}
& \|f_n - f^{[K_n]}\|_\infty  \, = \, \left\| \sum_{k=1}^{K_n} \left( \frac 1{\sqrt{2\pi n}} \alpha_{\theta(f),k} - \langle f,\varphi_k\rangle_{L_2} \right) \varphi_k \right\|_\infty
\\
& \, \leq \, \sqrt{\frac 2 \pi} \sum_{k=1}^{K_n} \left| \frac{\|M_k^*\|_F}{\sqrt{2\pi n}} - 1\right| \cdot | \langle f,\varphi_k \rangle_{L_2}| \, \leq \, \sqrt{\frac 2 \pi} \|f\|_{L_2} \cdot \left(\sum_{k=1}^{K_n} \left( 1- \sqrt{\frac{n-j'}{n}}  \right)^2  \right)^{1/2} \\
& \, \leq \, \sqrt{\frac 2 \pi} \frac{\sqrt{2 \pi}}{\rho^*} \cdot \left(\sum_{k=1}^{K_n} \frac{(j')^2}{ n^2}  \right)^{1/2} \leq \frac{\kappa_2 K_n^{1/2}}{n} =  {\cal O} (K_n/n) \,,
\end{align*}
for all $f \in W(s,L,\rho^*)$. Then we conclude using Lemma \ref{L:matrixf}(a) that
$$
\sup_{f \in W(s,L,\rho^*)} \|f - f_n\|_\infty \, = \, {\cal O}\big(K_n^{(1-s)/2} + K_n /n\big)\,.
$$
Moreover, $\hat f$ is a noisy version of $f_n$ and the noise has bounded support on the Euclidean sphere of radius $\gamma_n$. This gives
$$ \|f_n - \hat{f}\|_\infty \, \leq \, \frac1{\pi\sqrt{n}} \, \sum_{k=1}^{K_n} \big|\hat \alpha_k-\alpha_{\theta(f),k}\big| \, \leq \, \sqrt{K_n} \gamma_n / (\pi\sqrt{n})\,, $$
holds true by (\ref{eq:gamman}). Therefore, by (\ref{eq:Theo1}), when choosing $\gamma_n$ slightly larger than $K_n$ such that (\ref{eq:cond.pilot1}) and (\ref{eq:cond4'}) are satisfied, we may fix that $0 < \rho^*/2 \leq f_n,\hat{f} \leq 2/\rho^*$ for $n$ sufficiently large (uniformly with respect to $f$). Thus boundedness of the ratio $f_n/\hat{f}$ is guaranteed.

Let us now prove \eqref{equiv1} using the bounds on $f$ and $f_n$. We can bound
\begin{align*}
    & \|\log f - \log (f_n)\|_{L_2}^2  = \int \left( \int_{f_n(t,x)}^{f(t,x)} \frac 1u du \right)^2 \, dt \, dx \leq \int \left(\frac{(f-f_n)^2(t,x)}{f(t,x) f_n(t,x)} \right)^2 \, dt \, dx \\
    &\leq \frac{(\rho^*)^4}4 \frac9{(\rho^*)^2} \|f-f_n \|_{L_2}^2  \leq \frac{9 (\rho^*)^2}{8\pi n}  \bigg(  \|a - \alpha_{\theta(f)}\|^2 + \sum_{k>K_n} a_k^2\bigg)
    \leq \mathcal{O}(1) \bigg(\frac{K_n^2}{n^2} + K_n^{-s}\bigg) ,
\end{align*}
using the notation and results from the proof of Lemma~\ref{lem:pilot}. Thus, \eqref{equiv1} is satisfied as soon as
\begin{equation} \label{equiv2}
K_n^2 /n \to 0 \text{ and }n K_n^{-s} \to 0.
\end{equation}

We consider the experiment ${\cal B}^*_n$, in which one observes a function $\hat{f}$ on the domain $[0,1]\times(-\pi,\pi]$ and -- conditionally on that -- a Gaussian process $Y(t,x)$, $t\in [0,1]$, $x\in (-\pi,\pi]$ that is driven by
$$
dY(t,x) \, = \, \frac{\sqrt{n}}{2\sqrt{\pi}} \cdot f_n(t,x) / \hat{f}(t,x) \, dt \, dx \, + \, dW(t,x)\,.
$$
We show that
$$
\frac n{4\pi} \| \log (f_n / \hat f) - ( f_n / \hat f - 1)\|_{L_2}^2 \to 0\,.
$$
We use the inequality $|\log(1+u) - u|\leq 2 u^2$ for all $u: |u|< 1/2$, to get for $n$ large enough:
\begin{align*}
    &\frac n{4\pi} \| \log (1 + (f_n / \hat f -1) ) - ( f_n / \hat f -1)\|_{L_2}^2
    \leq \frac n{4\pi} \| 2 ( f_n / \hat f -1)^2\|_{L_2}^2 \\
    &\leq \frac n{4\pi} \int \frac 8{(\rho^*)^2} \cdot  (\hat f - f_n )^4  \leq \mathcal{O}(1) \cdot n \int \left(\sum_{j=1}^{K_n} \frac 1{\sqrt{2 \pi n} } \tilde \eta_j \varphi_j(t,x) \right)^4 dt dx\\
    & \leq \mathcal{O}(1) \cdot \frac 1n \, \|\tilde \eta\|^2 \cdot (\int \sum_{j=1}^{K_n} \varphi_j^2)^2 \leq  \mathcal{O}(1) \cdot \frac {\gamma_n^2 K_n^2}n,
\end{align*}
which tends to 0 under the condition that $\gamma_n^2 K_n^2/n \to 0$. This is implied by our final stronger condition and together with \eqref{equiv1} proves that the experiments $(\mathcal{A}_n^*)'$ and $\mathcal{B}_n^*$ are equivalent.

By the Cameron-Martin theorem, the conditional image measure $P_{f_n,Y|\hat{f}}$ of $Y$ given $\hat{f}$ is dominated by the Wiener measure $P_W$ of the process $W$, and has the density
\begin{align*}
& p_{f_n,Y|\hat{f}}(W|\hat{f})  \, = \, \frac{dP_{f_n,Y|\hat{f}}}{dP_W}(W|\hat{f}) \\ & \, = \, \exp\big( - \frac n{4 \pi}\, \frac 12\, \big\|f_n/\hat{f}\big\|_2^2\big) \cdot \exp\Big(\frac1{2 \pi \sqrt{2 } } \, \sum_{j=1}^{K_n} \alpha_{\theta(f),j} \int_0^1\int_{-\pi}^\pi \varphi_j(t,x) / \hat{f}(t,x) \, dW(t,x)\Big)\,,
\end{align*}
so that, by the Fisher-Neyman factorization lemma, $(\hat{f},Y_{K_n})$ with
$$ Y_{K_n} \, := \, \Big\{\int_0^1\int_{-\pi}^\pi \varphi_j(t,x) / \hat{f}(t,x) \, dY(t,x)\Big\}_{j=1,\ldots,K_n}\,, $$
forms a sufficient statistic for $\theta=\theta(f)$. Therefore, ${\cal B}^*_n$ is equivalent to the experiment ${\cal C}^*_n$ in which only the statistic $(\hat{f},Y_{K_n})$ from the experiment ${\cal B}^*_n$ is observed. Since the $j$th component of $Y_{K_n}$ equals
\begin{align*}
Y_{K_n,j} \, = \, & \frac1{2 \pi \sqrt{2}} \sum_{j'=1}^{K_n} \alpha_{\theta(f),j'} \int_0^1\int_{-\pi}^\pi \varphi_j(t,x) \varphi_{j'}(t,x) / \hat{f}^2(t,x) \, dt \, dx \\ & \, + \, \int_0^1\int_{-\pi}^\pi \varphi_j(t,x)  / \hat{f}(t,x) \, dW(t,x)\,,
\end{align*}
the conditional distribution of $Y_{K_n}$ given $\hat{f}$ turns out to be
$$ {\cal N}\big( \Gamma_{\hat{f}} \alpha_\theta / (2 \pi \sqrt{2}), \Gamma_{\hat{f}}\big)\,, $$
with the matrix
$$ \Gamma_{\hat{f}} \, := \, \Big\{\int_0^1 \int_{-\pi}^\pi \varphi_j(t,x) \varphi_{j'}(t,x) / \hat{f}^2(t,x) \, dt \, dx\Big\}_{j,j'=1,\ldots,K_n}\,, $$
and the vector $\alpha_{\theta} = (\alpha_{\theta(f),j})_{j=1,\ldots,K_n}$. The matrix $\Gamma_{\hat{f}}$ is positive definite as
$$ {\bf v}^\top \Gamma_{\hat{f}} {\bf v} \, = \, \int_0^1\int_{-\pi}^\pi \Big|\sum_{j=1}^{K_n} v_j \varphi_j(t,x)\Big|^2 / \hat{f}^2(t,x) \, dt\, dx \, \geq \, \|{\bf v}\|^2 / \|\hat{f}\|_\infty^2\,, $$
for all $v=(v_j)_{j=1,\ldots,K_n}$. This inequality also implies that $$\|\Gamma_{\hat{f}}^{-1/2}\|_{sp} \leq \|\hat{f}\|_\infty.$$ Thus ${\cal C}_n^*$ is equivalent to the experiment ${\cal D}_n^*$ in which we observe $\hat{f}$ and -- conditionally on that -- a ${\cal N}(\alpha_\theta / (2 \pi \sqrt{2})  , \Gamma_{\hat{f}}^{-1})$-distributed random variable.
In the next step we need the proxy
$$ \widetilde{f^{-1/2}} \, := \, \sum_{j=1}^{K_n} \langle \hat{f}^{-1/2} , \varphi_j\rangle_{L_2} \cdot \varphi_j\,, $$
for $\hat{f}^{-1/2}$ to define the symmetric and positive semi-definite matrix
$$ \tilde{\Gamma}_{\hat{f}} \, := \,  \Big\{\int_0^1 \int_{-\pi}^\pi \varphi_j(t,x) \varphi_{j'}(t,x) \big\{\widetilde{f^{-1/2}}(t,x)\big\}^4 \, dt \, dx\Big\}_{j,j'=1,\ldots,K_n}\,. $$
\begin{lem} \label{L:f12}
For any $2 < s^* < s-1$ it holds that
$$ \sup_{f\in W(s,L,\rho^*)} \, \big\|\widetilde{f^{-1/2}} - \hat{f}^{-1/2}\big\|_\infty \, = \, {\cal O}\big(K_n^{1-s^*/2}\big)\,, $$
whenever
$\lim_{n\to\infty} \, K_n^{(s^*+1)/2} \cdot \gamma_n \cdot n^{-1/2} \, = \, 0$,
with $\gamma_n$ as in (\ref{eq:gamman}).
\end{lem} 
We have
\begin{align} \nonumber
\big\|&\Gamma_{\hat{f}}^{-1/2}(\Gamma_{\hat{f}} - \tilde{\Gamma}_{\hat{f}})\Gamma_{\hat{f}}^{-1/2}\big\|_F^2  \, \leq \, \|\Gamma_{\hat{f}}^{-1/2}\|_{sp}^4 \,\big\|\Gamma_{\hat{f}} - \tilde{\Gamma}_{\hat{f}}\|_F^2 \\ \nonumber
& \, \leq \, \|\hat{f}\|_\infty^4\, \sum_{j,j'=1}^{K_n} \big\langle \varphi_{j'} \, , \, \varphi_j  \big\{\big(\widetilde{f^{-1/2}}\big)^4 - (\hat{f}^{-1/2})^4 \big\}\big\rangle_{L_2}^2 \\ \nonumber
& \, \leq \, 16 (\rho^*)^{-4}\cdot K_n \cdot  \big\| \widetilde{f^{-1/2}}^4 - (\hat{f}^{-1/2})^4 \big\|_\infty^2 \\  \label{eq:sp.10}
& \, \leq \, 256 (\rho^*)^{-4} \cdot K_n \cdot \big\{\sqrt{2/\rho^*} +  \big\| \widetilde{f^{-1/2}} - \hat{f}^{-1/2} \big\|_\infty\big\}^6 \cdot  \big\| \widetilde{f^{-1/2}} - \hat{f}^{-1/2} \big\|_\infty^2\,.
\end{align}
Thus, under the conditions of Lemma \ref{L:f12}, the above term converges to zero uniformly with respect to $f$ with the rate $K_n^{3-s^*}$. Then we may also fix that $$\|\widetilde{f^{-1/2}}\|_\infty \leq \sqrt{3/\rho^*}$$ for $n$ sufficiently large (again uniformly with respect to $f$).
Combining Lemma \ref{L:sp.1} with (A.4) in \cite{Reiss2011}, we deduce that
\begin{align*}
& \mathbb{E}_\theta  \, H^2\big({\cal N}(\alpha_\theta/(2\pi \sqrt{2}) , \Gamma_{\hat{f}}^{-1}) \, , \, {\cal N}(\alpha_\theta/(2\pi \sqrt{2}) , \tilde{\Gamma}_{\hat{f}}^{-1})\big) \\ & \, \leq \, 2 \, \mathbb{E}_\theta \min\big\{1 , \big\|\Gamma_{\hat{f}}^{1/2}(\Gamma_{\hat{f}}^{-1} - \tilde{\Gamma}_{\hat{f}}^{-1})\Gamma_{\hat{f}}^{1/2}\big\|_F^2\big\}  \, = \, 2\, \int_0^1 \mathbb{P}_\theta\big[\big\|\Gamma_{\hat{f}}^{1/2}(\Gamma_{\hat{f}}^{-1} - \tilde{\Gamma}_{\hat{f}}^{-1})\Gamma_{\hat{f}}^{1/2}\big\|_F^2 > t\big] \, dt \\
& \, \leq \, 2\mathbb{P}_\theta\big[\big\|\Gamma_{\hat{f}}^{-1/2}(\Gamma_{\hat{f}} - \tilde{\Gamma}_{\hat{f}})\Gamma_{\hat{f}}^{-1/2}\big\|_{sp}>1/2\big] \\ & \quad + \, 2\, \int_0^1 \mathbb{P}_\theta\big[\Gamma_{\hat{f}}^{-1/2}(\Gamma_{\hat{f}} - \tilde{\Gamma}_{\hat{f}})\Gamma_{\hat{f}}^{-1/2}\big\|_{sp}\leq 1/2 \, , \,  \big\|\Gamma_{\hat{f}}^{1/2}(\Gamma_{\hat{f}}^{-1} - \tilde{\Gamma}_{\hat{f}}^{-1})\Gamma_{\hat{f}}^{1/2}\big\|_F^2 > t\big] \, dt \\
& \, \leq \, 2\mathbb{P}_\theta\big[\big\|\Gamma_{\hat{f}}^{-1/2}(\Gamma_{\hat{f}} - \tilde{\Gamma}_{\hat{f}})\Gamma_{\hat{f}}^{-1/2}\big\|_F>1/2\big] \\ & \hspace{5.3cm} + \, 2\, \int_0^1 \mathbb{P}_\theta\big[ 4 \big\| \Gamma_{\hat{f}}^{-1/2}(\Gamma_{\hat{f}} - \tilde{\Gamma}_{\hat{f}})\Gamma_{\hat{f}}^{-1/2}\big\|_F^2 > t\big] \, dt \\
& \, \leq \, 16 \, \mathbb{E}_\theta \, \big\|\Gamma_{\hat{f}}^{-1/2}(\Gamma_{\hat{f}} - \tilde{\Gamma}_{\hat{f}})\Gamma_{\hat{f}}^{-1/2}\big\|_F^2\,,
\end{align*}
where $H$ denotes the Hellinger distance. Thus, by (\ref{eq:sp.10}), asymptotic equivalence of ${\cal D}^*_n$ and the experiment ${\cal E}^*_n$, in which we observe $\hat{f}$ and -- conditionally on that -- a ${\cal N}( \alpha_\theta/(2\pi \sqrt{2}),\tilde{\Gamma}_{\hat{f}}^{-1})$-distributed random variable, follows whenever
\begin{equation} \label{eq:cond.sp.1}
\lim_{n\to\infty} \, K_n^{3-s^*} \, + \, K_n^{s^*+1} \cdot \gamma_n^2 \cdot n^{-1}  \, = \, 0\,.
\end{equation}

In this case, the inequality (\ref{eq:sp.10}) and the arguments from (\ref{eq:cond1.1}) yield that all eigenvalues of $\tilde{\Gamma}_{\hat{f}}$ have the lower bound $c / \|\hat{f}\|_\infty^2 \geq 4c (\rho^*)^{-2}$ for any fixed $c\in (0,1)$ when $n$ is sufficiently large. We have $\check{M}_j\cdot \sqrt{n/(2\pi)} = \Psi^{-1}(\varphi_j) \in {\cal M}(\kappa_1,\kappa_2)$ and define
$$
\widetilde{C^{-1/2}} := \Psi^{-1}(\widetilde{f^{-1/2}})\in {\cal M}(\kappa_1,\kappa_2)
$$ (where $\Psi=\Psi_{\kappa_1,\kappa_2}$). The $(j,j')$th entry of $\tilde{\Gamma}_{\hat{f}}$ can be written as
\begin{align*} & \, \big\langle \Psi(\widetilde{C^{-1/2}}) \, \Psi(\check{M}_j) \, \Psi(\widetilde{C^{-1/2}}) \, , \, \Psi(\widetilde{C^{-1/2}}) \, \Psi(\check{M}_{j'}) \, \Psi(\widetilde{C^{-1/2}}) \big\rangle_{{L_2^n}}   =   \, \big\langle \widetilde{f^{-1/2}}^4 \, \varphi_{j} \, , \, \varphi_{j'}\big\rangle_{L_2}\,.
\end{align*}
We introduce the $K_n\times K_n$-matrix $\check{\Gamma}_C$ whose $(j,j')$th entry equals
\begin{align*}
& \big\langle \widetilde{C^{-1/2}} \check{M}_j \widetilde{C^{-1/2}}  ,  \widetilde{C^{-1/2}} \check{M}_{j'} \widetilde{C^{-1/2}} \big\rangle_{F} = \big\langle \Psi(\widetilde{C^{-1/2}} \check{M}_j \widetilde{C^{-1/2}})  ,  \Psi(\widetilde{C^{-1/2}} \check{M}_{j'} \widetilde{C^{-1/2}}) \big\rangle_{L_2^n} \\
& = \big\langle \widetilde{f^{-1/2}}^4 \varphi_j , \varphi_{j'} \big\rangle_{L_2} + \big\langle \Psi(\widetilde{C^{-1/2}})  \Psi(\check{M}_j)  \Psi(\widetilde{C^{-1/2}}) , \delta_{j'} \big\rangle_{L_2^n} \\ & \quad + \big\langle \delta_{j} , \Psi(\widetilde{C^{-1/2}}) \Psi(\check{M}_{j'}) \Psi(\widetilde{C^{-1/2}})  \big\rangle_{L_2^n}
+ \langle \delta_j , \delta_{j'}\rangle_{L_2^n}\,,
  \end{align*}
as the map $\Psi=\Psi_{3\kappa_1,3\kappa_2}$ is an isometry, where
$$ \delta_j \, := \, \Psi(\widetilde{C^{-1/2}} \check{M}_j \widetilde{C^{-1/2}}) \, - \, \Psi(\widetilde{C^{-1/2}}) \, \Psi(\check{M}_j) \, \Psi(\widetilde{C^{-1/2}})\,. $$
It holds that
\begin{align*}
\big\|\tilde{\Gamma}_{\hat{f}}^{-1/2}(\tilde{\Gamma}_{\hat{f}} - \check{\Gamma}_C)\tilde{\Gamma}_{\hat{f}}^{-1/2}\big\|_F^2 & \, \leq \, \|\tilde{\Gamma}_{\hat{f}}^{-1/2}\|_{sp}^4 \, \big\|\tilde{\Gamma}_{\hat{f}} - \check{\Gamma}_C\big\|_F^2 \\ & \, \leq \, 3 \|{\hat{f}}\|_\infty^4 \cdot \Big\{ 2 \sum_{j=1}^{K_n} \big\|\widetilde{f^{-1/2}}^2 \delta_{j}\big\|_{L_2^n}^2 + \Big(\sum_{j=1}^{K_n} \|\delta_j\|_{L_2^n}^2\Big)^2\Big\} \\
& \, \leq \, 48 (\rho^*)^{-4} \cdot \Big\{ 18(\rho^*)^{-2} \sum_{j=1}^{K_n} \|\delta_{j}\|_{L_2^n}^2 + \Big(\sum_{j=1}^{K_n} \|\delta_j\|_{L_2^n}^2\Big)^2\Big\}\,, \end{align*}
when $n$ is sufficiently large (uniformly with respect to $f$). Note that
\begin{align*}
& \|\delta_j\|_{L_2^n} \, \leq  \,  \|\widetilde{f^{-1/2}}\|_\infty \cdot \big\| \Psi(\check{M}_j) \, \Psi(\widetilde{C^{-1/2}}) \, - \, \Psi(\check{M}_j \widetilde{C^{-1/2}})\big\|_{L_2^n} \\
&\ \ \ \ \  \, +\, \big\| \Psi(\widetilde{C^{-1/2}} \check{M}_j \widetilde{C^{-1/2}}) \, - \, \Psi(\widetilde{C^{-1/2}}) \Psi(\check{M}_j \widetilde{C^{-1/2}})\big\|_{ L_2^n } \\
& \, \leq \, 2\pi \sqrt{3/\rho^*} \cdot \|\check{M}_j\|_F\, \|\widetilde{C^{-1/2}}\|_F \cdot \kappa_1\kappa_2 \cdot n^{-3/2}  +  2\pi \|\widetilde{C^{-1/2}}\|_F \, \|\check{M}_j \widetilde{C^{-1/2}}\|_F \, \kappa_1\kappa_2  n^{-3/2}
\end{align*}
by Lemma A.1 where $\|\check{M}_j\|_F = 1$, $\|\widetilde{C^{-1/2}}\|_F = \|\widetilde{f^{-1/2}}\|_{L_2^n} \leq \sqrt{3n/\rho^*}$ and
\begin{align*} \big\|\check{M}_j \widetilde{C^{-1/2}}\big\|_F & \, = \, \big\|\Psi\big(\check{M}_j \widetilde{C^{-1/2}}\big)\big\|_{L_2^n} \\ & \, \leq \, \big\|\Psi(\check{M}_j)\cdot \Psi\big(\widetilde{C^{-1/2}}\big)\big\|_{L_2^n} \, + \, 2\pi \, \big\|\widetilde{C^{-1/2}}\big\|_F \, \kappa_1\kappa_2 \, n^{-3/2} \\
& \, = \, \big\|\widetilde{f^{-1/2}}\big\|_\infty + 2\pi\, \|\widetilde{f^{-1/2}}\|_{L_2^n}\, \kappa_1\kappa_2\, n^{-3/2} \, \leq \, \sqrt{3/\rho^*} \cdot \big(1 + 2\pi \kappa_1\kappa_2 / n\big)\,,
\end{align*}
so that $\|\delta_j\|_{L_2^n}^2$ is bounded from above by some global constant factor times $K_n^2/n^2 + K_n^4/n^4$. Again using Lemma \ref{L:sp.1} we conclude asymptotic equivalence of ${\cal E}^*_n$ and the experiment ${\cal F}^*_n$ in which we observe $\hat{f}$ and -- conditionally -- on that a ${\cal N}(\alpha_\theta/(2 \pi \sqrt{2}), \check{\Gamma}_C^{-1})$-distributed random variable whenever
\begin{align} \label{eq:cond.sp.2}
\lim_{n\to\infty}\, K_n^3 n^{-2} + K_n^5 n^{-4} \,  = \, 0,
\end{align}
which is satisfied under (\ref{eq:Theo1}). The matrices $\check{M}_j$ and, hence, $\widetilde{C^{-1/2}}$ are symmetric and real-valued. The $(j,j')$th entry of $\check{\Gamma}_C$ may be written as
\begin{align*}
\mbox{tr}&\big(\widetilde{C^{-1/2}}^2 \check{M}_j \widetilde{C^{-1/2}}^2 \check{M}_{j'}\big) \, = \, \big\langle |\widetilde{C^{-1/2}}| \check{M}_j |\widetilde{C^{-1/2}}| \, , \,  |\widetilde{C^{-1/2}}| \check{M}_j |\widetilde{C^{-1/2}}|\big\rangle_F\,, \end{align*}
where $|\widetilde{C^{-1/2}}|$ stands for the positive square root of the matrix $\widetilde{C^{-1/2}}^2$. We use the analogous arguments from Section \ref{GOE} in order to show equivalence to the experiment ${\cal G}_n^*$ in which we observe $\hat{f}$ and -- conditionally on that --
$$ -|\widetilde{C^{-1/2}}/\sqrt{a_n^*}| {\check{C}}_\theta |\widetilde{C^{-1/2}}/\sqrt{a_n^*}| \, + \, G\,, $$
with the GOE $G$, ${\check{C}}_\theta = \sum_{j=1}^{K_n} \alpha_{\theta(f),j} \check{M}_j$ and $a_n^*:= a_n \sqrt{\pi n}$. Introducing the matrix $\check{C} := \sum_{j=1}^{K_n} \hat \alpha_{j} \check{M}_j$, which is fully known by observing $\hat{f}$ (and so is $\widetilde{C^{-1/2}}$), the model ${\cal G}_n^*$ is equivalent to the experiment ${\cal H}_n^*$ in which we observe $\hat{f}$ and -- conditionally on $\hat{f}$ --
$$  |\widetilde{C^{-1/2}}/\sqrt{a_n^*}| (\check{C}-{\check{C}}_\theta) |\widetilde{C^{-1/2}}/\sqrt{a_n^*}| \, + \, G\,. $$
This proves the following theorem.
\begin{theo} \label{AE_BWN}
     Let $\mathcal{A}^*_n$ denote the bivariate Gaussian white noise model described by
    $$
    dX(t,x) = \log f(t,x) \, dt \, dx + 2\sqrt{\frac{\pi}{n}} \, dW(t,x), \quad \text{on } [0,1]\times(-\pi,\pi],
    $$
    for $f$ belonging to the class $W(s,L,\rho^*)$ with $s>10$ and let the circulant-type matrices $\{ \check M_k; \, k=1,...,K_n\}$ be defined in Appendix A of the supplementary material. If $K_n$ is such that $K_n^{10} \log n/n \to 0$ as $n\to \infty$, then the experiments  $\mathcal{A}^*_n$ and $\mathcal{H}^*_n$ are equivalent.
\end{theo}

\subsection{Connecting the GOE Models}

In this section we will show asymptotic equivalence between the experiment ${\cal K}_n$ from Section \ref{GOE}, see also Theorem~\ref{AE_LSP}, and the experiment ${\cal H}_n^*$ from Section \ref{LSDM}. By the Bretagnolle-Huber inequality and Jensen's inequality, the Le Cam distance between the experiment ${\cal K}_n$ and ${\cal L}_n$ has the upper bound
\begin{align*}
\Big\{1 - \exp\Big(- \mathbb{E}_\theta\,  \big\| |\widetilde{C^{-1/2}}/\sqrt{a_n^*}| \, \check{\Delta}_\theta \, |\widetilde{C^{-1/2}}/\sqrt{a_n^*}| - C^{-1/2} \Delta_\theta C^{-1/2}\big\|_F^2/4\Big)\Big\}^{1/2}\,,
\end{align*}
where $\check{\Delta}_\theta := \check{C} - {\check{C}}_\theta$ since the Kullback-Leibler distance between the Gaussian Orthogonal Ensemble with the shift $|\widetilde{C^{-1/2}}/\sqrt{a_n^*}| \, \check{\Delta}_\theta \, |\widetilde{C^{-1/2}}/\sqrt{a_n^*}|$, on the one hand, and $C^{-1/2} \Delta_\theta C^{-1/2}$, on the other hand, equals
$$ \big\||\widetilde{C^{-1/2}}/\sqrt{a_n^*}|\, \check{\Delta}_\theta \, |\widetilde{C^{-1/2}}/\sqrt{a_n^*}| - C^{-1/2} \Delta_\theta C^{-1/2}\big\|_F^2/4\,. $$
Consider that
\begin{align} \nonumber
& \big\||\widetilde{C^{-1/2}}/\sqrt{a_n^*}| \, \check{\Delta}_\theta \, |\widetilde{C^{-1/2}}/\sqrt{a_n^*}| - C^{-1/2} \Delta_\theta C^{-1/2}\big\|_F^2 \\  \nonumber & \, = \, \big\||\widetilde{C^{-1/2}}/\sqrt{a_n^*}| \, \check{\Delta}_\theta \,  |\widetilde{C^{-1/2}}/\sqrt{a_n^*}| - C^{-1/2} \Delta_\theta C^{-1/2}\big\|_F^2 \\ & \nonumber
 \leq 3 (a_n^*)^{-1} \cdot \big\||\widetilde{C^{-1/2}}/\sqrt{a_n^*}| \, - \, C^{-1/2}\big\|_F^2\cdot  \|\check{\Delta}_\theta\|_{sp}^2\cdot \|\widetilde{C^{-1/2}}\|_{sp}^2 \\ \nonumber & \quad + 3 (a_n^*)^{-1} \cdot \|C^{-1/2}\|_{sp}^2 \cdot \big\| \check{\Delta}_\theta - \Delta_\theta\big\|_F^2 \cdot \|\widetilde{C^{-1/2}}\|_{sp}^2  \\ \label{eq:GOEbound}& \quad +  3 \cdot \| C^{-1/2}\|_{sp}^2\cdot \|\Delta_\theta\|_{sp}^2\cdot \big\||\widetilde{C^{-1/2}}/\sqrt{a_n^*}| - C^{-1/2}\big\|_F^2,
\end{align}
where we have used that $\| |\widetilde{C^{-1/2}}|\|_{sp} = \|\widetilde{C^{-1/2}}\|_{sp}$. As we have already seen in Section \ref{DR} it holds that $\|C^{-1/2}\|_{sp} \leq 1/\sqrt{c\rho}$ for any $c\in (0,1)$ and $n$ sufficiently large (uniformly in $f$). Moreover,
\begin{align*}
& \sup_{f\in W(s,L,\rho^*)} \, \|\widetilde{C^{-1/2}}\|_{sp}^2 \, \leq \, \sup_{f\in W(s,L,\rho^*)}\, \Big(\sum_{j=1}^{K_n} |\langle \hat{f}^{-1/2} , \varphi_j\rangle_{L_2}| \cdot \|\check{M}_j\|_{sp}\cdot \sqrt{n/(2\pi)}\Big)^2 \\ &  \, \leq \, \frac{n}{2\pi} \, \sup_{f\in W(s,L,\rho^*)}  \big\|\hat{f}^{-1/2}\big\|_{L_2}^2 \, \sum_{j=1}^{K_n} \big(2\|\check{M}_j - M_j\|_F^2 + 2 \|M_j\|_{sp}^2\big) \, = \, {\cal O}\big(K_n^{3/2}\big)\,,
\end{align*}
by the Cauchy-Schwarz inequality, (\ref{eq:MMcheck}) and the fact that $\|\check{M}_j\|_{sp} \leq 2\sqrt{2/\pi}$. Furthermore,
\begin{align*}
\big\|&|\widetilde{C^{-1/2}}/\sqrt{a_n^*}| - C^{-1/2}\big\|_F^2  \, \leq \, \|C^{-1/2}\|_{sp}^2 \, \big\|I_n - |\widetilde{C^{-1/2}}/\sqrt{a_n^*}| C^{1/2}\big\|_F^2 \\ & \, \leq \,  \|C^{-1/2}\|_{sp}^2 \, \big\|I_n - \widetilde{C^{-1/2}}^2 C/a_n^*\big\|_F^2  \\
& \, \leq \, 2 \, \|C^{-1/2}\|_{sp}^2 \, \big\|I_n - \widetilde{C^{-1/2}}^2 \check{C}/a_n^*\big\|_F^2 \, + \, 2(a_n^*)^{-2} \, \|C^{-1/2}\|_{sp}^2 \,  \|\widetilde{C^{-1/2}}\|_{sp}^4  \|C - \check{C}\|_F^2,
\end{align*}
where
\begin{align*}
\big\|&I_n - \widetilde{C^{-1/2}}^2 \check{C} / a_n^*\big\|_F  \, = \, \big\|1 - \Psi\big(\widetilde{C^{-1/2}}^2 \check{C}\big)/a_n^*\big\|_{L_2^n}  \\ &\, \leq \, \big\|1 - 2\pi \cdot\widetilde{f^{-1/2}}^2\cdot \hat{f}/a_n^*\big\|_{L_2^n} \, + \, \big\|\Psi^2\big(\widetilde{C^{-1/2}}\big) \Psi(\check{C}) - \Psi\big(\widetilde{C^{-1/2}}^2\big) \Psi(\check{C})\big)\big\|_{L_2^n}/a_n^*  \\ & \quad\  + \, \big\|\Psi\big(\widetilde{C^{-1/2}}^2\big) \Psi(\check{C}) - \Psi\big(\widetilde{C^{-1/2}}^2 \check{C}\big)\big\|_{L_2^n}/a_n^*\,,
\end{align*}
as $\Psi(\check{C}) = 2\pi \hat{f}$. Now we fix that $a_n^*=2\pi$ (thus $a_n = 2\sqrt{\pi} / \sqrt{n}$). By Lemma A.1 we deduce that
\begin{align*}
& \sup_{f\in W(s,L,\rho^*)}\, \big\|I_n - \widetilde{C^{-1/2}}^2 \check{C} / a_n^*\big\|_F^2 \, \leq \, {\cal O}\big(K_n^2 n^{-1} + K_n^4 n^{-3}\big) \\ &  \, + \, \sup_{f\in W(s,L,\rho^*)} \, 2n \cdot \|\hat{f}\|_\infty^2 \cdot \big\|\hat{f}^{-1/2} - \widetilde{f^{-1/2}}\big\|_\infty^4   + 4n \cdot \|\hat{f}\|_\infty^2 \cdot \|\hat{f}^{-1/2}\|_\infty^2   \big\|\hat{f}^{-1/2} - \widetilde{f^{-1/2}}\big\|_\infty^2 \\
& \, = \, {\cal O}\big(K_n^2 n^{-1} + K_n^4 n^{-3} + n K_n^{2-s^*}\big)\,,
\end{align*}
where we apply Lemma \ref{L:f12} under the given constraints. Moreover, using the Cauchy-Schwarz inequality and (\ref{eq:MMcheck}), we find that
\begin{align*}
\sup_{f\in W(s,L,\rho^*)}\, 2(a_n^*)^{-2} \, \|C^{-1/2}\|_{sp}^2 \, \|\widetilde{C^{-1/2}}\|_{sp}^4  \|C - \check{C}\|_F^2 \, = \, {\cal O}\big(K_n^{9/2}\big)\,.
 \end{align*}
It follows from (\ref{eq:gamman}) and (\ref{eq:condMk}) that
$$ \sup_{f\in W(s,L,\rho^*)} \big\|\Delta_{\theta(f)}\big\|_{sp}^2 \, = \, {\cal O}\big(\gamma_n^2 K_n n^{-1}\big)\,. $$
Using (\ref{eq:MMcheck}) we derive that
\begin{align*}
\sup_{f\in W(s,L,\rho^*)} \|\Delta_{\theta(f)} - \check{\Delta}_{\theta(f)}\|_F^2 & \, = \, \sup_{f\in W(s,L,\rho^*)} \Big\|\sum_{j=1}^{K_n} (\alpha_{\theta(f),j} - \alpha_j) \, (\check{M}_j - M_j)\Big\|_F^2 \\
& \, \leq \, \gamma_n^2 \cdot \sum_{j=1}^{K_n} \|\check{M}_j - M_j\|_F^2 \, = \, {\cal O}\big(\gamma_n^2 K_n^{3/2} n^{-1}\big)\,.
\end{align*}
From there it follows that
$$ \sup_{f\in W(s,L,\rho^*)} \big\|\check{\Delta}_{\theta(f)}\big\|_{sp}^2 \, = \, {\cal O}\big(\gamma_n^2 K_n^{3/2} n^{-1}\big)\,. $$
Piecing together these results, we can show that the term (\ref{eq:GOEbound}) converges to zero uniformly with respect to $f\in W(s,L,\rho^*)$ whenever (\ref{eq:Theo1}) is granted and $s>10$. This choice of $s$ enables us to impose that $s^*=7$ in Lemma \ref{L:f12} and the condition therein is granted under our assumptions. 
Thus, together with Theorems~\ref{AE_LSP} and~\ref{AE_BWN}, this proves the following theorem.

\begin{theo} \label{AE}
   The experiments $\mathcal{A}_n$ of the locally stationary process with spectral density $f$ and $\mathcal{A}^*_n$ of the bivariate Gaussian white noise with drift function $\log f$ are asymptotically equivalent for $f$ belonging to the class $W(s,L,\rho^*)$ with $s>10$, $L>0$ and $\rho^*$ in (0,1).
\end{theo} 

%
%





The supplemental material of \cite{AELS_main} is organized as follows. In Section \ref{Circ}, the circulant-type matrices are introduced and properties of the linear map $\Psi$ from the circulant-type matrices to $L_2^n$ are deduced. Section \ref{sec:proofs} contains the proofs of the lemmata in the main article as well as the proof of Lemma \ref{L:almosthom}.

\begin{acks}[Acknowledgments]
The authors acknowledge the CIRM, France, for sponsoring this research project under the grant 3212/2024. We are grateful to the reviewers for their very constructive comments.
\end{acks}
\begin{funding}
A. Meister and A. Rohde acknowledge support of the Research Unit 5381 (DFG), ME 2114/5-1 and RO 3766/8-1, respectively.
\end{funding}

\newpage

\appendix
\section{Circulant-type Matrices $\{\check{M}_k\}_k$} \label{Circ}

As a companion matrix for circulant matrices we specify the permutation matrix
$$ \check{M} \, := \, \begin{pmatrix} 0 & 1 & 0 & \cdots & 0 \\
                              \vdots &  &  & \ddots & \vdots \\
                              0 & 0 & 0 & \cdots & 1 \\
                              1 & 0 & 0 & \cdots & 0 \end{pmatrix}\,. $$
It is well known that $\check{M}$ has the eigenvalues $\lambda(j) := \exp(2\pi i j / n)$ and the corresponding eigenvectors
$$ {\bf v}_j \, := \, \big(\lambda(0),\lambda(j),\ldots,\lambda(j(n-1))\big)^\top\,. $$
Let $\Lambda$ be the diagonal matrix which contains the components of ${\bf v}_1$ as the entries of its main diagonal. Moreover define the matrices
$$ \check{M}_{0,j,j'} \, := \, \Lambda^j \check{M}^{j'}\,, \quad j,j'=0,\ldots,n-1\,. $$
Note the $\check{M}_{0,j,j'}$ are also defined for all integer $j,j'$ and are $n$-periodic in both indices $j,j'$ (and so are the $\lambda(j)$ and ${\bf v}_j$ in their index). We deduce that
\begin{align*}
\check{M}_{0,j,j'} {\bf v}_k & \, = \, \lambda(j'k) \cdot \Lambda^j {\bf v}_k \, = \, \lambda(j'k) \cdot\big(\lambda(0),\lambda(j+k),\ldots,\lambda((j+k)(n-1))\big)^\top \\
& \, = \, \lambda(j'k) \cdot {\bf v}_{j+k}\,. \end{align*}
For $j,j',k,k'=-\lfloor (n-1)/2\rfloor,\ldots,\lfloor(n-1)/2\rfloor$, it holds that
\begin{align*}
\langle \check{M}_{0,j,j'} , \check{M}_{0,k,k'} \rangle_F & \, = \, \mbox{tr}\big(\Lambda^j \check{M}^{j'} \big[\Lambda^{-k} \check{M}^{k'}\big]^\top\big) \, = \, \mbox{tr}\big(\Lambda^j \check{M}^{j'} \check{M}^{-k'} \Lambda^{-k}\big) \, = \, \mbox{tr}\big(\Lambda^{j-k} \check{M}^{j'-k'}\big)\,.
\end{align*}
If $j'\neq k'$ then all entries of the main diagonal of $\check{M}^{j'-k'}$ vanish; and so do the main diagonal entries of $\Lambda^{j-k} \check{M}^{j'-k'}$ so that $\check{M}_{0,j,j'}$ and $\check{M}_{0,k,k'}$ are orthogonal. Otherwise, for $j'=k'$, note that
\begin{align*}
 \mbox{tr}\big(\Lambda^{j-k} \check{M}^{j'-k'}\big) & \, = \, \mbox{tr}\big(\Lambda^{j-k}\big) \, = \, \sum_{\ell=0}^{n-1} \exp(2\pi i (j-k)/  n)^\ell \, = \, \begin{cases} n\,, & \mbox{ for }j=k\,, \\
0\,,&  \mbox{ for }j\neq k\,. \end{cases}
\end{align*}
Therefore the $\check{M}_{0,j,j'}$, $j,j'=-\lfloor(n-1)/2\rfloor,\ldots,\lfloor(n-1)/2\rfloor$, form an orthogonal system with respect to the Frobenius inner product in the space of all complex $n\times n$-matrices and satisfy $\|\check{M}_{0,j,j'}\|_F^2=n$ for all $j,j'$. Moreover we deduce that
\begin{align} \nonumber
\check{M}_{0,j_1,j_1'} \check{M}_{0,j_2,j_2'} {\bf v}_k & \, = \, \lambda(j_2' k) \, \check{M}_{0,j_1,j_1'} {\bf v}_{j_2+k} \, = \, \lambda(j_2'k)\cdot \lambda(j_1'(j_2+k)) \, {\bf v}_{j_2+k+j_1} \\ \nonumber
& \, = \, \lambda(j_1'j_2) \cdot \lambda((j_2'+j_1')k) {\bf v}_{k+j_1+j_2} \, = \, \lambda(j_1'j_2) \cdot \check{M}_{0,j_1+j_2,j_1'+j_2'} {\bf v}_k\,,
\end{align}
for all $k=0,\ldots,n-1$. As the ${\bf v}_k$, $k=0,\ldots,n-1$, form a basis of $\mathbb{C}^n$ it holds that
\begin{equation} \label{eq:almosthom} \check{M}_{0,j_1,j_1'} \check{M}_{0,j_2,j_2'} \, = \, \lambda(j_1'j_2) \cdot \check{M}_{0,j_1+j_2,j_1'+j_2'}\,. \end{equation}
The linear hull of the $\check{M}_{0,j,j'}$, $j=-\kappa_{1},\ldots,\kappa_{1}$, $j' = -\kappa_{2},\ldots,\kappa_{2}$ with $\kappa_\ell < \lfloor(n-1)/2\rfloor$ for $\ell=1,2$ is denoted by ${\cal M}(\kappa_{1},\kappa_{2})$. Moreover, recall we defined the functions $\tilde{\varphi}_{j,j'}$ by
$$\tilde{\varphi}_{j,j'}(t,x) \, := \, \exp\big(2\pi i j t + i j' x\big)\,. $$
They form an orthogonal basis in the Hilbert space $L_2^n$ of all measurable, complex-valued and squared integrable functions on the domain $[0,1]\times (-\pi,\pi]$ which has the inner product
$$ \langle f,g \rangle_{{L_2^n}} \, := \, \frac{n}{2\pi} \int_0^1\int_{-\pi}^\pi f(t,x) \, \overline{g(t,x)} \, dt\, dx \, = \, \frac{n}{2\pi} \, \langle f,g\rangle_{L_2}\,, $$
and they satisfy $\|\tilde{\varphi}_{j,j'}\|_{L_2^n}^2 = n$ for all integer $j,j'$. There exists a unique linear map  $\Psi_{\kappa_1,\kappa_2}$ on the domain ${\cal M}(\kappa_1,\kappa_2)$ which maps $\check{M}_{0,j,j'}$ to $\tilde{\varphi}_{j,j'}$ for all $j=-\kappa_{1},\ldots,\kappa_{1}$, $j' = -\kappa_{2},\ldots,\kappa_{2}$. When the codomain of $\Psi_{\kappa_1,\kappa_2}$ is restricted to its range it forms an isometric isomorphism between two Hilbert spaces. On the other hand, $\Psi_{2\kappa_1,2\kappa_2}$ is also an approximate homomorphism with respect to standard matrix multiplication in ${\cal M}(\kappa_1,\kappa_2)$ and pointwise multiplication of functions in $L_2^n$. That is made concrete in the following lemma. 
\begin{lem}\label{L:almosthom}
Set $\Psi=\Psi_{2\kappa_1,2\kappa_2}$ for $\kappa_\ell < \lfloor (n-1)/2\rfloor / 2$ and $\ell=1,2$. For all $A,B \in {\cal M}(\kappa_1,\kappa_2)$, it holds that
$$ \|\Psi(AB) - \Psi(A)\Psi(B)\|_{L_2^n}^2 \, \leq \, 4 \pi^2 \, \|A\|_F^2 \|B\|_F^2 \cdot \kappa_1^{2} \kappa_2^{2} / n^3 \,.$$
\end{lem}
We define $\check{M}_{j,j'}^\pm := \Psi^{-1}\big(\varphi_{j,j'}^\pm\big)$, $j\leq \kappa_1$, $j'\leq \kappa_2$, with $\Psi$ as in Lemma \ref{L:almosthom} and $\varphi_{j,j'}^\pm$ as in (\ref{eq:varphi+}) and (\ref{eq:varphi-}). Note that the $\varphi_{j,j'}^\pm$ are located in the range of the linear map $\Psi$ in view of equations \eqref{eq:tildevarphi+} and \eqref{eq:tildevarphi-}. 
Since $$\check{M}_{0,-k,k'} = \overline{\check{M}_{0,k,k'}}\ \ \text{and}\ \  \check{M}_{0,k,-k'} = \check{M}_{0,k,k'}^\top$$ the matrices $\check{M}_{k,k'}^\pm$ are symmetric and real-valued. As $\Psi$ is isometric it holds that
$$ \langle \check{M}_{j,j}^\pm , \check{M}_{k,k'}^\pm\rangle_F \, = \, \langle \varphi_{j,j'}^\pm , \varphi_{k,k'}^\pm\rangle_{L_2^n} \, = \, \frac{n}{2\pi} \cdot  \langle \varphi_{j,j'}^\pm , \varphi_{k,k'}^\pm\rangle_{L_2}\,, $$
so that the matrices $\check{M}_{j,j'}^\pm$ form an orthogonal system and satisfy $\big\|\check{M}_{j,j'}^\pm\big\|_F^2 = n/(2\pi)$. We fix some enumeration $\check{M}_k$, $k=1,\ldots,K_n$, of the full collection of matrices $\check{M}_{j,j'}^+ \cdot \sqrt{2\pi/n}$, $j=0,\ldots,\kappa_1$, $j'=0,\ldots,\kappa_2$, and the $\check{M}_{j,j'}^- \cdot \sqrt{2\pi/n}$, $j=1,\ldots,\kappa_1$, $j'=0,\ldots,\kappa_2$, where $K_n = (2\kappa_1+1)(\kappa_2+1)$. Thus the $\check{M}_k$, $k=1,\ldots,K_n$, form an orthonormal system with respect to the Frobenius inner product.

\section{Proofs of the Lemmata} \label{sec:proofs}

\noindent {\it Proof of Lemma \ref{L:Gamma_theta}:} We start with showing the positive definiteness of $\Gamma$.
\begin{align} \nonumber
\inf \big\{t^\top \Gamma t \, : \, \|t\|=1\big\} & \, = \, 2 \, \inf \big\{\mbox{tr}\big(M^{[t]} C^{-1} M^{[t]} C^{-1}\big) \, : \, \|t\|=1\big\} \\ \nonumber
& \, = \, 2\, \inf \big\{\|C^{-1/2} M^{[t]} C^{-1/2}\|_F^2 \, : \, \|t\|=1\big\} \\ \nonumber
& \, = \, 2\, \inf \Big\{\sum_{j,j'=1}^n \lambda_j^{-1} \lambda_{j'}^{-1} \, \big(\varphi_j^\top  M^{[t]} \varphi_{j'}\big)^2 \, : \, \|t\|=1\Big\} \\ \nonumber
& \, \geq \, 2\, \|C\|_{sp}^{-2} \cdot \inf \big\{ \|M^{[t]}\|_F^2 \, : \, \|t\|=1\big\} \\ \nonumber
& \, = \, 2\, \|C\|_{sp}^{-2} \cdot \inf \Big\{ \sum_{k=1}^{K_n} t_k^2 \|M_k\|_F^2 \, : \,  \|t\|=1\Big\} \\  \nonumber
& \, = \, 2\, \|C\|_{sp}^{-2} \\ \label{eq:Proof1}
& \, \geq \, c \rho^2\,,
\end{align}
for $n$ sufficiently large, where $M^{[t]} := \sum_{k=1}^{K_n} t_k \cdot M_k$; and the $\varphi_j$ and the $\lambda_j$ denote the eigenvectors and the eigenvalues of $C$, respectively. This yields positive definiteness (and, hence, invertibility) of the symmetric matrix $\Gamma$. We deduce that
\begin{align} \nonumber
\big\|\Gamma^{-1/2} (\Gamma_\theta - \Gamma) \Gamma^{-1/2}\big\|_F   \, \leq \, \|C\|_{sp}^{2} \cdot &\big[ 2 \, \big\|\big\{\mbox{tr}\big(\Delta_\theta C^{-1} M_{k'} C^{-1} M_k C^{-1}\big)\big\}_{k,k'} \big\|_F \\ \nonumber  & \, + \, \big\|\big\{\mbox{tr}\big(\Delta_\theta C^{-1} M_k C^{-1} \Delta_\theta C^{-1} M_{k'} C^{-1}\big)\big\}_{k,k'}\big\|_F\big]\,.
\end{align}
Moreover consider that
\begin{align*}
\big\| \big\{\mbox{tr}\big(\Delta_\theta C^{-1} M_{k'} C^{-1} M_k C^{-1}\big)\big\}_{k,k'} \big\|_F^2 & \, = \, \big\| \big\{\big\langle M_k, C^{-1} M_{k'} C^{-1} \Delta_\theta C^{-1}\big\rangle_F \big\}_{k,k'}\big\|_F^2 \\
& \, \leq \, \sum_{k'=1}^{K_n} \big\|C^{-1} M_{k'} C^{-1} \Delta_\theta C^{-1}\big\|_F^2 \\
& \, \leq \, \sum_{k'=1}^{K_n} \|C^{-1}\|_{sp}^6 \cdot \|\Delta_\theta\|_{sp}^2 \cdot \|M_{k'}\|_F^2 \\
& \, \leq \, K_n \cdot c^{-6} \cdot \rho^{-6} \cdot \gamma_n^2 \cdot \sum_{k=1}^{K_n} \|M_k\|_{sp}^2\,,
\end{align*}
for all $\theta$ and any fixed constant $c\in (0,1)$ for $n$ sufficiently large, where we use (\ref{eq:gamman}) and $\|M_k\|_F=1$ for all $k=1,\ldots,K_n$.   Furthermore, the same arguments provide that, for all $\theta\in\Theta$,
\begin{align*}
\big\| \big\{\mbox{tr}\big(\Delta_\theta C^{-1} M_k C^{-1} \Delta_\theta C^{-1} M_{k'} C^{-1}\big)\big\}_{k,k'}\big\|_F^2 & \, \leq \, \sum_{k'=1}^{K_n} \| C^{-1} \Delta_\theta C^{-1} M_{k'} C^{-1} \Delta_\theta C^{-1}\|_F^2 \\
& \, \leq \, K_n \cdot c^{-8} \cdot \rho^{-8} \cdot \gamma_n^4 \cdot \Big(\sum_{k=1}^{K_n} \|M_k\|_{sp}^2\Big)^2 \,,
\end{align*}
again for any fixed constant $c\in (0,1)$ and $n$ sufficiently large. Analogously we derive that
\begin{align} \nonumber
\big\|\Gamma^{-1/2} (\tilde{\Gamma}_\theta - \Gamma) \Gamma^{-1/2}\big\|_F  \, \leq \, & \|C\|_{sp}^2 \cdot\big[  2 \big\| \big\{\mbox{tr}\big(M_k C^{-1} \Delta_\theta C_\theta^{-1} M_{k'} C_\theta^{-1}\big)\big\}_{k,k'}\big\|_F \\ \nonumber
& \, + \, \big\|\big\{\mbox{tr}\big(M_k C^{-1} \Delta_\theta C_\theta^{-1} M_{k'} C^{-1} \Delta_\theta C_\theta^{-1} C^{-1}\big)\big\}_{k,k'}\big\|_F\big],
\end{align}
and
\begin{align*}
 \big\| \big\{\mbox{tr}\big(M_k C^{-1} \Delta_\theta C_\theta^{-1} M_{k'} C_\theta^{-1}\big)\big\}_{k,k'}\big\|_F^2 & \, \leq \, \sum_{k=1}^{K_n} \big\|C_\theta^{-1} M_k C^{-1} \Delta_\theta C_\theta^{-1}\big\|_F^2 \\
& \, \leq \, \sum_{k=1}^{K_n} \|C_\theta^{-1}\|_{sp}^4 \cdot \|C^{-1}\|_{sp}^2 \cdot\|\Delta_\theta\|_{sp}^2 \cdot \|M_k\|_F^2 \\
& \, \leq \, K_n \cdot c^{-6} \cdot \rho^{-6} \cdot \gamma_n^2 \cdot \sum_{k=1}^{K_n} \|M_k\|_{sp}^2\,,
\end{align*}
as well as
\begin{align*}
\big\|\big\{\mbox{tr}\big(M_k C^{-1} \Delta_\theta & C_\theta^{-1} M_{k'} C^{-1} \Delta_\theta C_\theta^{-1} C^{-1}\big)\big\}_{k,k'}\big\|_F^2 \\ &  \leq  \sum_{k'=1}^{K_n} \big\|C^{-1} \Delta_\theta C_\theta^{-1} M_{k'} C^{-1} \Delta_\theta C_\theta^{-1} C^{-1}\big\|_F^2 \\
& \, \leq \, K_n\cdot c^{-8} \cdot \rho^{-8} \cdot \gamma_n^4 \cdot \Big(\sum_{k=1}^{K_n} \|M_k\|_{sp}^2\Big)^2 \,.
\end{align*}
Using that $\|C\|_{sp}^2 \leq (c\rho)^{-2}$ for any $c\in (0,1)$ when $n$ is sufficiently large, we can prove that
$$ \lim_{n\to\infty} \, \sup_\theta \big\{\big\|\Gamma^{-1/2} (\Gamma-\Gamma_\theta) \Gamma^{-1/2}\big\|_F \, + \, \big\|\Gamma^{-1/2} (\Gamma-\tilde{\Gamma}_\theta) \Gamma^{-1/2}\big\|_F\big\} \, = \, 0\,, $$
under the conditions of the lemma. The second claim, i.e. the lower bound on the eigenvalues, follows from (\ref{eq:Proof1}) and the first claim using the arguments from (\ref{eq:cond1.1}). 
\hfill $\square$ \\[0.5cm]

\noindent{\it Proof of Lemma \ref{L:FouTails}}: In order to derive an upper bound on the tails of $\Psi_n$ we consider, by (\ref{eq:Gauss.1}), that
$$ |\Psi_n(t_1,\ldots,t_{K_n})| \, = \, \prod_{j=1}^n \big(1 + 4 \lambda_{j,\theta,t}^2\big)^{-1/4}\,, $$
for all $t \in \mathbb{R}^{K_n}$. Since
$$ \sum_{j=1}^n \lambda_{j,\theta,t}^2 \, = \, \big\|D_\theta^{[t]}\big\|_F^2 \, = \, \frac12 t^\top \Gamma_\theta t\,, $$
and
$$ \lambda_{j,\theta,t}^2 \, \leq \, \big\|D_\theta^{[t]}\big\|_{sp}^2\,, $$
for all $j=1,\ldots,n$, we deduce
$$ \# \Big\{j=1,\ldots,n \, : \, \lambda_{j,\theta,t}^2 \geq \frac14 t^\top \Gamma_\theta t / n\Big\} \, \geq \, \frac14 t^\top \Gamma_\theta t / \big\|D_\theta^{[t]}\big\|_{sp}^2\,, $$
so that
$$ |\Psi_n(t_1,\ldots,t_{K_n})| \, \leq \, \big(1 + t^\top \Gamma_\theta t / n\big)^{-t^\top \Gamma_\theta t / (16 \|D_\theta^{[t]}\|_{sp}^2)}\,, $$
for all $t \in \mathbb{R}^{K_n}$ with $D_\theta^{[t]}\neq 0$ and, hence,
$$ |\Psi_n^*(t_1,\ldots,t_{K_n})| \, \leq \, \big(1 + \|t\|^2 / n\big)^{- \|t\|^2 / (16 \|\tilde{D}_\theta^{[t]}\|_{sp}^2)}\,, $$
for all $t \in \mathbb{R}^{K_n}$ with $\tilde{D}_\theta^{[t]}\neq 0$, thus, for all $t\in \mathbb{R}^{K_n} \backslash\{0\}$ since the matrices $\sqrt{2} \tilde{D}_{\theta,k}$, $k=1,\ldots,K_n$, form an orthonormal system with respect to the Frobenius inner product. Using (\ref{eq:bound_Dt}) it follows that
$$ \|t\|^2 / \big(16 \|\tilde{D}_\theta^{[t]}\|_{sp}^2\big) \, \geq \, 1/ (16\mu_{n}^2)\,, $$
for all $t\neq 0$. So the following upper bound applies,
$$ |\Psi_n^*(t_1,\ldots,t_{K_n})| \, \leq \, \big(1 + \|t\|^2 / n\big)^{- 1/(16 \mu_{n}^2)}\,, \qquad \forall t\in \mathbb{R}^{K_n}\,. $$
For any $R\geq 0$ we consider the integral
\begin{align} \nonumber
& \idotsint_{\|t\|\geq R} \, |\Psi_n^*(t_1,\ldots,t_{K_n})| \, dt_1\cdots dt_{K_n} \\  \nonumber
& \, \leq \, n^{K_n/2}  \, \idotsint_{\|t\|\geq R n^{-1/2}} \, \big(1 + \|t\|^2\big)^{- 1/(16\mu_{n}^2)}\, dt_1\cdots dt_{K_n} \\  \nonumber
& \, = \, 2 \, (n \pi)^{K_n/2} \, [\Gamma(K_n/2)]^{-1} \, \int_{r \geq R n^{-1/2}} \, \big(1 + r^2\big)^{-  1/(16\mu_{n}^2)} \, r^{K_n-1}\, dr \\  \label{eq:Gauss10}
& \, = \,  (n \pi)^{K_n/2} \, [\Gamma(K_n/2)]^{-1} \, \int_{r \geq R^2 / n} \, \big(1 + r\big)^{-  1/(16\mu_{n}^2)} \, r^{K_n/2-1}\, dr\,.
\end{align}
Putting $R=0$, the above term is finite whenever $\mu_n^{-2} > 8 K_n$.  If even $\mu_n^{-2} > 8K_n+16$ it follows that (\ref{eq:Gauss10}) is bounded from above as claimed in the lemma. \hfill $\square$ \\[0.5cm]

\noindent {\it Proof of Lemma \ref{L:Edge}}: Expanding the $\exp$-function by its power series and using some appropriate factorization we deduce that
\begin{align*}
\exp&\Big\{\frac12 \sum_{\ell=3}^\infty \frac1\ell (2i)^\ell \cdot \mbox{tr}\big[\big(\tilde{D}_\theta^{[t]}\big)^\ell\big]\Big\}  \\ & \, = \, 1 \, + \, \sum_{k=1}^\infty \frac1{2^k k!} \,\Big( \sum_{\ell=3}^\infty \frac1\ell (2i)^\ell \sum_{j_1,\ldots,j_\ell=1}^{K_n} t_{j_1} \cdots t_{j_\ell} \cdot \mbox{tr}\big(\tilde{D}_{\theta,j_1} \cdots \tilde{D}_{\theta,j_\ell})\Big)^k \\
& \, = \, 1 \, + \, \sum_{q=3}^\infty (2i)^q \sum_{j_1,\ldots,j_q=1}^{K_n} t_{j_1} \cdots t_{j_q} \cdot \sum_{k=1}^{\lfloor q/3\rfloor} \frac1{2^k k!}\, \sum_{\ell_1,\ldots,\ell_k=3}^\infty {\bf 1}_{\{q\}}(\ell_1+\cdots+\ell_k) \\ & \hspace{5.5cm} \cdot\prod_{m=0}^{k-1} \frac1{\ell_{m+1}} \cdot \mbox{tr}\big(\tilde{D}_{\theta,j_{L(m)+1}} \cdots \tilde{D}_{\theta,j_{L(m+1)}}\big)\,,
\end{align*}
where $L(m) := \sum_{j=1}^m \ell_j$ whenever $\|t\|< 1/(2\mu_n)$. Thus, for these $t$, the remainder term admits the following multivariate Taylor expansion
$$ \exp\Big\{\frac12 \sum_{\ell=3}^\infty \frac1\ell (2i)^\ell \cdot \mbox{tr}\big[\big(\tilde{D}_\theta^{[t]}\big)^\ell\big]\Big\} \, = \, 1 \, + \, \sum_{q=3}^\infty \sum_{\|{\bf m}\|_1=q} \nu_{\bf m} \cdot t_1^{m_1} \cdots t_{K_n}^{m_{K_n}}\,, $$
with ${\bf m} := (m_1,\ldots,m_{K_n}) \in \mathbb{N}_0^{K_n}$; $\|{\bf m}\|_1 = m_1 +\cdots+m_{K_n}$; and the coefficients
\begin{align*} \nu_{\bf m} & \, := \, (2i)^q \, \sum_{j_1,\ldots,j_q=1}^{K_n} \Big\{\prod_{j=1}^{K_n} {\bf 1}_{\{m_j\}}\big(\#\{p : j_p=j\}\big)\Big\} \\ & \cdot \sum_{k=1}^{\lfloor q/3\rfloor} \frac1{2^k k!}\, \sum_{\ell_1,\ldots,\ell_k=3}^\infty {\bf 1}_{\{q\}}(\ell_1+\cdots+\ell_k) \cdot\prod_{m=0}^{k-1} \frac1{\ell_{m+1}} \cdot \mbox{tr}\big(\tilde{D}_{\theta,j_{L(m)+1}} \cdots \tilde{D}_{\theta,j_{L(m+1)}}\big)\,, \end{align*}
when $q = \|{\bf m}\|_1$. Using (\ref{eq:bound_Dt}) for $t={\bf e}_k$, i.e. the unit vector with $1$ as its $k$th component, we find that
\begin{align*} \|\tilde{D}_{\theta,k}\|_{sp} & \, \leq \, \|C_\theta^{1/2} C^{-1}\|_{sp} \, \|C^{-1} C_\theta^{1/2}\|_{sp} \, \sum_{k'=1}^{K_n} \big|\big(\Gamma_\theta^{-1/2}\big)_{k,k'}\big|\, \|M_{k'}\|_{sp}
 \, \leq \, \mu_n\,,
\end{align*}
for all $k=1,\ldots,K_n$. Then we derive the following upper bound on the involved traces,
\begin{align*}
\big|&\mbox{tr}\big(\tilde{D}_{\theta,j_{L(m)+1}} \cdots \tilde{D}_{\theta,j_{L(m+1)}}\big)\big| \, = \, \big| \big\langle \tilde{D}_{\theta,j_{L(m)+1}} \cdots \tilde{D}_{\theta,j_{L(m+1)-1}} \, , \, \tilde{D}_{\theta,j_{L(m+1)}} \big\rangle_F\big| \\
& \, \leq \, \big\|\tilde{D}_{\theta,j_{L(m)+1}} \cdots \tilde{D}_{\theta,j_{L(m+1)-1}}\big\|_F \cdot \big\| \tilde{D}_{\theta,j_{L(m+1)}}\big\|_F \\
& \, \leq \, \big\|\tilde{D}_{\theta,j_{L(m)+1}}\big\|_{sp} \cdots \big\|\tilde{D}_{\theta,j_{L(m+1)-2}}\big\|_{sp} \cdot \big\|\tilde{D}_{\theta,j_{L(m+1)-1}}\big\|_F \cdot  \big\|\tilde{D}_{\theta,j_{L(m+1)}}\big\|_F \\
& \, \leq \, \frac12 \, \mu_n^{\ell_{m+1}-2}\,,
\end{align*}
so that
\begin{align*}
\prod_{m=0}^{k-1} \big|\mbox{tr}\big(\tilde{D}_{\theta,j_{L(m)+1}} \cdots \tilde{D}_{\theta,j_{L(m+1)}}\big)\big| & \, \leq \, 2^{-k} \, \mu_n^{q-2k}\,.
\end{align*}
Moreover it holds that
\begin{align*}
\sum_{\ell_1,\ldots,\ell_k=3}^\infty & {\bf 1}_{\{q\}}(\ell_1+\cdots+\ell_k) \cdot\prod_{m=0}^{k-1} \frac1{\ell_{m+1}} \, \leq \, \Big(\sum_{\ell=3}^q 1/\ell\Big)^k \, \leq \, (\log q)^k\,.
\end{align*}
Piecing together these inequalities the following bound on the coefficients applies,
\begin{align} \nonumber
|\nu_{\bf m}| & \, \leq \, 2^q \, \sum_{j_1,\ldots,j_q=1}^{K_n} \Big\{\prod_{j=1}^{K_n} {\bf 1}_{\{m_j\}}\big(\#\{p : j_p=j\}\big)\Big\} \cdot \sum_{k=1}^{\lfloor q/3\rfloor} \frac1{4^k k!} \, (\log q)^k \, \mu_n^{q-2k} \\ \nonumber
& \, \leq \, 2^q \, \mu_n^{q/3} \, q^{1/4} \, \sum_{j_1,\ldots,j_q=1}^{K_n} \Big\{\prod_{j=1}^{K_n} {\bf 1}_{\{m_j\}}\big(\#\{p : j_p=j\}\big)\Big\} \\ \nonumber
& \, = \,  {q \choose m_1,\ldots,m_{K_n}} \, 2^q \cdot \mu_n^{q/3} \cdot q^{1/4}\,,
\end{align}
when $\|{\bf m}\|_1=q$, using the multinomial theorem in the last step. Finally the residual term $R_{Q,n}(t)$ can be bounded as follows,
\begin{align} \nonumber
|&R_{Q,n}(t)|  \, \leq \, \sum_{q>Q} \sum_{\|{\bf m}\|_1=q} |\nu_{\bf m}| \cdot |t_1|^{m_1} \cdots |t_{K_n}|^{m_{K_n}} \\  \nonumber
& \, \leq \, \sum_{q>Q}  2^q \cdot \mu_n^{q/3} \cdot q^{1/4} \cdot \|t\|_1^q \, \sum_{\|{\bf m}\|_1=q} {q \choose m_1,\ldots,m_{K_n}} \, (|t_1|/\|t\|_1)^{m_1} \cdots (|t_{K_n}|/\|t\|_1)^{m_{K_n}} \\  \nonumber
& \, \leq \,  \sum_{q>Q}  \big(2 K_n^{1/2} \mu_n^{1/3} \cdot \|t\|\big)^q \cdot q^{1/4}  \\  \nonumber
& \, \leq \, \sum_{q>Q} \big\{2 K_n^{1/2} \mu_n^{1/3} (Q+1)^{1/(4Q+4)} \cdot \|t\| \big\}^q \\ \nonumber
& \, = \, (Q+1)^{1/4} \cdot \big(2 K_n^{1/2} \mu_n^{1/3} \cdot \|t\| \big)^{Q+1} / \big\{1 - 2 K_n^{1/2} \mu_n^{1/3} (Q+1)^{1/(4Q+4)} \cdot \|t\|\big\}\,,
\end{align}
whenever $Q\geq 2$ and $\|t\| < (Q+1)^{1/(4Q+4)} \cdot\mu_n^{-1/3} K_n^{-1/2} /2$. Therein, $\|t\|_1 := \sum_{k=1}^{K_n} |t_k|$, and we have used the multinomial theorem again. \hfill $\square$ \\[0.5cm]


\noindent{\it Proof of Lemma~\ref{L:piecing}}: 
We use the decomposition of the Fourier integral to see that the probability in (\ref{eq:TVmin}) is bounded from below by
\begin{align} \nonumber
\mathbb{P}&\big[\tau^{[1]}_{n,\theta}(X) \, > \, (s-1) \, (2\pi)^{-K_n/2} \exp\big(-\|X\|^2/2\big) \, + \, \tau^{[2]}_{n,\theta} \, + \, \tau^{[3]}_{n,\theta} \, + \, \tau^{[4]}_{n,\theta}\big] \\  \nonumber
& \, \geq \, \mathbb{P}\big[\tau^{[1]}_{n,\theta}(X) \, > \, (s-1) \, (2\pi)^{-K_n/2} \exp\big(-\|X\|^2/2\big) / 2 \, , \\  \nonumber & \hspace{1.3cm} \, \tau^{[2]}_{n,\theta} \, + \, \tau^{[3]}_{n,\theta} \, + \, \tau^{[4]}_{n,\theta} \, < \, (1-s) \, (2\pi)^{-K_n/2} \exp\big(-\|X\|^2/2\big) / 2\big] \\  \nonumber
& \, \geq \, 1 \, - \, \mathbb{P}\big[\tau^{[1]}_{n,\theta}(X) \, \leq \, (s-1) \, (2\pi)^{-K_n/2} \exp\big(-\|X\|^2/2\big) / 2\big] \\ \label{eq:TV10} & \hspace{1cm} \, - \, \mathbb{P}\big[ \tau^{[2]}_{n,\theta} \, + \, \tau^{[3]}_{n,\theta} \, + \, \tau^{[4]}_{n,\theta} \, \geq \, (1-s) \, (2\pi)^{-K_n/2} \exp\big(-\|X\|^2/2\big) / 2\big]\,.
\end{align}
With the focus on the first term in (\ref{eq:TV10}) consider that
\begin{align*} \tau_{n,\theta}^{[1]}(x) & \, = \, \mbox{Re}\, (2\pi)^{-K_n/2}\sum_{q=3}^{Q_n} i^q \sum_{\|{\bf m}\|_1=q} \nu_{\bf m} \, \frac{\partial^q}{\partial x_1^{m_1} \cdots \partial x_{K_n}^{m_{K_n}}}\,  \exp\big(-\|x\|^2/2\big) \\
& \, = \, (2\pi)^{-K_n/2}  \exp\big(-\|x\|^2/2\big) \cdot \mbox{Re}\, \sum_{q=3}^{Q_n} (-i)^q \sum_{\|{\bf m}\|_1=q} \nu_{\bf m} \, \prod_{k=1}^{K_n} H_{m_k}(x_k)\,,
\end{align*}
where $H_j$ denotes the $j$th Hermite polynomial, so that
\begin{align} \nonumber
\mathbb{P}&\big[\tau^{[1]}_{n,\theta}(X) \, \leq \, (s-1) \, (2\pi)^{-K_n/2} \exp\big(-\|X\|^2/2\big) / 2\big] \\ \nonumber
& \, \leq \, \mathbb{P}\Big[\Big| \sum_{q=3}^{Q_n} (-i)^q \sum_{\|{\bf m}\|_1=q} \nu_{\bf m} \, \prod_{k=1}^{K_n} H_{m_k}(X_k)\Big| \, \geq \, (1-s)/2\Big] \\ \nonumber
& \, \leq \, 4 (1-s)^{-2} \cdot \mathbb{E} \Big| \sum_{q=3}^{Q_n} (-i)^q \sum_{\|{\bf m}\|_1=q} \nu_{\bf m} \, \prod_{k=1}^{K_n} H_{m_k}(X_k)\Big|^2 \\ \nonumber
& \, = \, 4(1-s)^{-2} \, \sum_{q,q'=3}^{Q_n} (-i)^{q-q'} \sum_{\|{\bf m}\|_1=q}  \sum_{\|{\bf m}'\|_1=q'} \nu_{\bf m} \cdot \overline{\nu_{{\bf m}'}} \cdot \prod_{k=1}^{K_n} \mathbb{E} \, H_{m_k}(X_k) H_{m'_k}(X_k) \\ \label{eq:TV20}
& \, = \, 4(1-s)^{-2} \, \sum_{q=3}^{Q_n} \sum_{\|{\bf m}\|_1=q} |\nu_{\bf m}|^2 \cdot \prod_{k=1}^{K_n} m_k!\,,
\end{align}
by Markov's inequality and the orthogonality of the Hermite polynomials. Then, the upper bound on the $\nu_{{\bf m}}$ from Lemma \ref{L:Edge} and the multinomial theorem yield that (\ref{eq:TV20}) is bounded from above by
\begin{align} \nonumber 4(&1-s)^{-2} \, \sum_{q=3}^{Q_n} (4 \mu_n^{2/3})^q \, q^{1/2}\, \sum_{\|{\bf m}\|_1=q} {q \choose m_1,\ldots,m_{K_n}}^2  \cdot \prod_{k=1}^{K_n} m_k! \\  \nonumber
& \, = \, 4(1-s)^{-2} \, \sum_{q=3}^{Q_n} q! \, (4 \mu_n^{2/3})^q \, q^{1/2}\, \sum_{\|{\bf m}\|_1=q} {q \choose m_1,\ldots,m_{K_n}} \\ \label{eq:TV30}
& \, = \, 4(1-s)^{-2} \, \sum_{q=3}^{Q_n} q! \, \big(4 K_n \mu_n^{2/3}\big)^q \, q^{1/2}\,,
\end{align}
which converges to $0$ whenever
\begin{equation*} \lim_{n\to\infty} \, K_n^3 \mu_n^2 \, = \, 0 \,, \end{equation*}
and
\begin{equation*} \label{eq:condCLT2} Q_n \, \leq \, \sqrt3 K_n^{-1} \mu_n^{-2/3} / 8  \,, \end{equation*}
as (\ref{eq:condCLT2}) guarantees that the summands in (\ref{eq:TV30}) decrease in $q=3,\ldots,{Q_n}$. These conditions, however, follow from (\ref{eq:cond.pilot1}), (\ref{eq:cond4'}) and (\ref{eq:Qn}) that we recall here:
\begin{equation*} 
\lim_{n\to\infty} \, Q_n \, K_n^{4/3} n^{-1/3}\, = \, 0\,.
\end{equation*}

The second term in (\ref{eq:TV10}) is bounded as follows,
\begin{align*}
\mathbb{P}&\big[(1-s) (2\pi)^{-K_n/2} \exp\big(-\|X\|^2/2\big) / 2 \, \leq \, \tau_{n,\theta}^{[2]} + \tau_{n,\theta}^{[3]} + \tau_{n,\theta}^{[4]}\big] \\
& \, \leq \, \mathbb{P}\Big[\sum_{k=1}^{K_n} X_k^2 \, \geq \,  -2 \Big(\log \frac6{1-s}\Big) \, - \, 2 \max_{j=2,3,4} \log\big\{(2\pi)^{K_n/2} \tau_{n,\theta}^{[j]}\big\}\Big]\,,
\end{align*}
so that, by Markov's inequality, this term converges to zero (uniformly w.r.t. $\theta$) whenever
\begin{equation*}
\max_{j=2,3,4} \sup_\theta \, \big( \log \tau_{n,\theta}^{[j]} \big) / K_n \, \longrightarrow \, - \infty\,, \end{equation*}
as $n\to\infty$.

Using Lemma \ref{L:FouTails} and putting $R=R_n$ we conclude that
\begin{align*} \big( \log \tau_{n,\theta}^{[4]} \big) / K_n & \, \leq \, - (\log 2) - \frac12 (\log \pi) - \{\log \Gamma(K_n/2)\}/K_n + \frac12 (\log n) \\ & \qquad - \Big(\frac1{16} K_n^{-1} \mu_n^{-2} - \frac12 - K_n^{-1}\Big) \cdot \log\big(1+R_n^2/n\big)\,, \end{align*}
which tends to $-\infty$ if
\begin{equation*} 
 R_n^2/n \, \longrightarrow \, 0\,,
\end{equation*}
and $R_n^2 K_n^{-1} n^{-1} \mu_n^{-2} \, / \, \log n \, \longrightarrow \, +\infty$ as $n\to\infty$, where the latter condition holds under
\begin{equation*} 
R_n^2 K_n^{-2}\, / \, \log n \, \longrightarrow \, +\infty\,.
\end{equation*}
Note that, under (\ref{eq:cond.pilot1}) and (\ref{eq:cond4'}), $R_n$ can always be chosen such that the two previous conditions, i.e. condition  (\ref{eq:condCLT4}) of the lemma, hold simultaneously.

Now we study the term $\tau_{n,\theta}^{[2]}$. We put $\nu_0:=1$ and $\nu_{\bf m} :=0$ for all $\|{\bf m}\|_1\in \{1,2\}$ and derive the upper bound
\begin{align*}
& \tau_{n,\theta}^{[2]}   \leq (2\pi)^{-K_n} \sum_{q=0}^{Q_n} \sum_{\|{\bf m}\|_1=q}| \nu_{\bf m}| \int_{\|t\|\geq R_n} \exp\big(-\|t\|^2/2\big) |t_1|^{m_1} \cdots |t_{K_n}|^{m_{K_n}}\, dt \\
& \leq  (2\pi)^{-K_n} \sum_{q=0,q\geq3}^{Q_n} \max\{1,q^{1/4}\} \, \big(2 \mu_n^{1/3} K_n^{1/2}\big)^q \, \int_{\|t\|\geq R_n} \exp\big(-\|t\|^2/2\big) \|t\|^{q} dt \\
& = 2 \pi^{K_n/2}  (2\pi)^{-K_n} \sum_{q=0,q\geq3}^{Q_n} \max\{1,q^{1/4}\} \big(2 \mu_n^{1/3} K_n^{1/2}\big)^q \\ & \hspace{6.8cm} \cdot \int_{r>R_n} \exp(-r^2/2) r^{q+K_n-1} dr / \Gamma(K_n/2) \\
& \leq \pi^{-K_n/2} \, R_n^{-1}  [\Gamma(K_n/2)]^{-1}  2^{1-K_n}  \\ & \hspace{3.6cm} \cdot \sum_{q=0,q\geq3}^{Q_n} \max\{1,q^{1/4}\} \big(2 \mu_n^{1/3} K_n^{1/2}\big)^q \sup_{r\geq R_n} \exp(-r^2/2) r^{q+K_n+1},
\end{align*}
when $R_n\geq 1$, using the upper bound on the $\nu_{\bf m}$ from Lemma \ref{L:Edge} and the multinomial theorem. Therein, the sum over $q$ is taken from $0$ to $Q_n$ where $q=1,2$ are omitted. Thus, under the constraint
\begin{equation*} 
R_n^2 \, \geq \, Q_n+K_n+1\,,
\end{equation*}
the supremum is taken at $r=R_n$ and the term $\tau_{n,\theta}^{[2]}$ admits the upper bound
\begin{align*}  2^{1-K_n} \pi^{-K_n/2} \,  [\Gamma(K_n/2)]^{-1} \, \sum_{q=0,q\geq3}^{Q_n} \max\{1,q^{1/4}\}\, \big(2 \mu_n^{1/3} K_n^{1/2} R_n\big)^q \, \exp\big(-R_n^2/2\big)\, R_n^{K_n}\,. \end{align*}
Thus $\sup_\theta \, (\log \tau_{n,\theta}^{[2]})/K_n$ diverges to $-\infty$ if (\ref{eq:condCLT7}) holds, that is: 
\begin{align*} 
& \frac12 R_n^2 K_n^{-1} \, - \, \log R_n \, \longrightarrow \, +\infty\,, 
\end{align*}
and if $\limsup_{n\to\infty} \, 2 \mu_n^{1/3} K_n^{1/2} R_n \, < \, 1\,$ as $n\to\infty$. The latter condition is satisfied when (\ref{eq:condCLT6}) holds, that is: 
\begin{equation*} 
 \lim_{n\to\infty} \, K_n^4 \, R_n^6 / n \, = \, 0\,.
\end{equation*}

Finally consider the term $\tau_{n,\theta}^{[3]}$. Note that
$$ \tau_{n,\theta}^{[3]} \, \leq \, (2\pi)^{-K_n/2} \cdot \sup_{\|t\|<R_n} |R_{Q,n}(t)|\,, $$
since $$ \idotsint \exp(-\|t\|^2/2) \, dt_1\cdots dt_{K_n} \, = \, (2\pi)^{K_n/2}\,. $$
If (\ref{eq:condCLT6}) and (\ref{eq:condCLT8})
are granted, then Lemma \ref{L:Edge} provides that
$$ \sup_\theta \, \big(\log \, \tau_{n,\theta}^{[3]}\big) / K_n \, \longrightarrow \, -\infty. $$
 \hfill $\square$ \\[0.5cm]

\noindent {\it Proof of Lemma \ref{L:sp.1}}: As $\|A^{-1/2}BA^{-1/2}-I\|_{sp} < 1$, where $I$ stands for the identity matrix of the corresponding size, the matrix $A^{-1/2}BA^{-1/2}$ is invertible (expandable by Neumann series) and so is $B$. Thus $B$ is positive definite. We deduce that
\begin{align*}
A^{1/2}&(A^{-1}-B^{-1})A^{1/2}  =  A^{-1/2}(B-A)B^{-1}A^{1/2}  =  A^{-1/2}(B-A)A^{-1/2} A^{1/2} B^{-1} A^{1/2} \\ & \, = \,  A^{-1/2}(B-A)A^{-1/2} \, + \,  A^{-1/2}(B-A)A^{-1/2} (A^{1/2} B^{-1} A^{1/2} - I) \\
& \, = \, A^{-1/2}(B-A)A^{-1/2} \, + \,  A^{-1/2}(B-A)A^{-1/2} \big(A^{1/2} (B^{-1} - A^{-1})A^{1/2}\big)\,.
\end{align*}
Applying the Frobenius norm we arrive at
\begin{align*}
\|&A^{1/2}(A^{-1}-B^{-1})A^{1/2}\|_F \\ & \, \leq \, \|A^{-1/2}(B-A)A^{-1/2}\|_F \, + \,  \|A^{-1/2}(B-A)A^{-1/2}\|_{sp}\, \|A^{1/2} (A^{-1} - B^{-1})A^{1/2}\|_F\,,
\end{align*}
which completes the proof of the lemma. \hfill $\square$ \\[0.5cm]

\noindent {\it Proof of Lemma \ref{lem: 3.1}}:
First note that as a consequence of the Cauchy-Schwarz inequality, $A,A_{t,n}^0\in S^2(s,\tilde L)$ for some $s>5/2$ implies that \begin{itemize}
\item $u\mapsto A(u,x)$ is Lipschitz continuous  with Lipschitz semi-norm uniformly bounded in $x\in[-\pi,\pi]$ and 
\item $x\mapsto f(u,x)$ as well as the products $x\mapsto A_{t,n}^0(x)\overline{A_{s,n}^0(x)}$ possess some Sobolev smoothness larger than $3/2$, uniformly in $u\in[0,1]$ and $(s,t)\in\{1,\dots, n\}^2$, respectively. 
\end{itemize}
We start with proving $
\Arrowvert \vartheta^n-\theta(f) \Arrowvert_F=o(1)$. By part (ii) of the definition,
\begin{align*}
\sum_{r=1}^n\sum_{t=1}^{n}(\theta_{r,t}-\vartheta_{r,t}^n)^2
&=2\sum_{t=1}^n\sum_{r=1}^{t-1}(\theta_{r,t}-\vartheta_{r,t}^n)^2+\mathcal{O}\Big(\frac{1}{n}\Big).
\end{align*}
With
\begin{align*}
I_1(r,t,n;x)&:=A\Big(\frac{(r\wedge t)-1}{n},x\Big)\overline{A\Big(\frac{(r\wedge t)-1}{n},x\Big)}-A(r/n,x)\overline{A(t/n,x)},\\
I_2(r,t,n;x)&:=\Big[A(r/n,x) - A_{r,n}^0(x)\Big]\overline{A^0_{t,n}}\ \ \text{and}\\ 
I_3(r,t,n;x)&:=A(r/n,x) \Big[ \overline{A(t/n,x)} - \overline{A_{t,n}^0(x)}\Big],
\end{align*}
we get
\begin{align*}
\theta_{r,t}&-\vartheta_{r,t}^n=\int_{-\pi}^\pi \exp\big(ix(r-t)\big)\big( I_1(r,t,n;x)+I_2(r,t,n;x)+I_3(r,t,n;x)\big)dx.
\end{align*}
By (ii) and the Lipschitz continuity of $A$ in its first coordinate, we find for any $1\leq m\leq n$ 
\begin{align*}
\sum_{t=1}^n\sum_{k=1}^{m\wedge (t-1)}&(\theta_{t-k,t}-\vartheta_{t-k,t}^n)^2 \\ & 
=\sum_{t=1}^n\sum_{k=1}^{m\wedge (t-1)}\bigg\{\int_{-\pi}^{\pi}\exp(-ixk)\bigg(\sum_{i=1}^3I_i(t-k,t,n;x)\bigg)dx\bigg\}^2\\
&\lesssim \Arrowvert f\Arrowvert_{\infty}\frac{m}{n}+\sum_{t=1}^n\sum_{k=1}^{m\wedge( t-1)}\bigg\{\int_{-\pi}^{\pi}\exp(-ixk)I_1(t-k,t,n;x)dx\bigg\}^2\\ 
&\lesssim \Arrowvert f\Arrowvert_{\infty}\frac{m}{n} +\sum_{t=1}^n\sum_{k=1}^{m\wedge (t-1)}\Arrowvert f\Arrowvert_{\infty}\frac{k^2}{n^2} \\
&\lesssim \Arrowvert f\Arrowvert_{\infty}\frac{m^3}{n},
\end{align*}
where $\lesssim$ means less or equal up to some multiplicative real constant which does not depend on the variable parameters in the expression. Next, with $a_{f,u}(k):=\int_{-\pi}^\pi f(u,x)\exp(ixk)dx$ and $a_{g}(k):=\int_{-\pi}^\pi g(x)\exp(ixk)dx$,  
\begin{align*}
\sum_{t=1}^n\sum_{k=m\wedge (t-1)}^{t-1}(\theta_{r,t}-\vartheta_{r,t}^n)^2&\leq 2\sum_{t>m}^n\sum_{k=m}^{t-1}\bigg\{\big\arrowvert a_{f, \frac{t-k-1}{n}}(k)\big\arrowvert^2+ \big\arrowvert a_{A_{t-k,n}^0}*a_{\overline{A_{t,n}^0}}(k)\big\arrowvert^2\bigg\}\\
&\lesssim  \frac{n}{m^{2\gamma}}
\end{align*}
 for some $\gamma>3/2$ by the Sobolev smoothness of $A_{s,n}^0(\cdot)\overline{A_{t,n}^0}(\cdot)$  and $f(u,\cdot)$ in the second coordinate. The first  bound yields the requirement $m=o(n^{1/3})$. Choosing $m=n^{\alpha}$ with $1/(2\gamma)<\alpha<1/3$  yields $\Arrowvert\theta(f)-\vartheta^n\Arrowvert_F=o(1)$. 

As this implies in particular $\Arrowvert\theta(f)-\vartheta^n\Arrowvert_{sp}=o(1)$, the statement on the extremal eigenvalues then follows together with Lemma 4.4 in \cite{Dahlhaus1996}.
$~$ \hfill $\square$ \\[0.5cm]

\noindent {\it Proof of Lemma \ref{L:matrixf}}: (a)\; We deduce that
\begin{align*}
\|f - f^{[K_n]}\|_\infty & \, \leq \, \sum_{k>K_n} |\langle f,\varphi_k\rangle_{L_2}| \cdot \|\varphi_k\|_\infty  \\
& \, \leq \, (2/\sqrt{\pi}) \, \sqrt{L} \,  \left(\sum_{j,j'} \, {\bf 1}\{j>\kappa_1\mbox{ or }j'>\kappa_2\} \, \big(j^2+(j')^2\big)^{-s} \right)^{1/2} \\
& \, \leq \, (2/\sqrt{\pi}) \, \sqrt{L} \, \left(\iint_{x^2+y^2>\min\{\kappa_1,\kappa_2\}^2} (x^2 + y^2)^{-s} \, dx \, dy \right)^{1/2}\\
& \, \leq \, (2/\sqrt{\pi}) \, \sqrt{L} \, \left( \int_0^{\pi/2} \int_{r>\min\{\kappa_1,\kappa_2\}} (r^2)^{-s} \, r \, dr \, d\theta \right)^{1/2}\\
& \, \leq \,  \sqrt{2L/(2-2s)} \,  \min\{\kappa_1,\kappa_2\} ^{1-s} \,, 
\end{align*}
for all $f\in W(s,L,\rho^*)$ with $s>1$. 

Concerning the $L_2$ norm, 
\begin{align*}
\|f - f^{[K_n]}\|^2_{L_2} & \, = \, \sum_{k>K_n} |\langle f,\varphi_k\rangle_{L_2}|^2  \\
& \, \leq \,  \sum_{j,j'} \, {\bf 1}\{j>\kappa_1\mbox{ or }j'>\kappa_2\} |\langle f,\varphi_{j,j'} \rangle_{L_2}|^2 \, \big(j^2+(j')^2\big)^{s} \cdot \min\{\kappa_1,\kappa_2\}^{-2s} \\
& \, \leq \,  {L} \, \min\{\kappa_1,\kappa_2\}^{-2s}\,, 
\end{align*}
for all $f\in W(s,L,\rho^*)$ with $s>1$. 

\noindent (b)\; If $f$ and $\tilde f$ belong to $W(s,L,\rho^*)$ with $s>1$, then $h : =f-\tilde f$ is also uniformly bounded by $2\pi / \rho^*$ and satisfies the ellipsoid condition. Thus,
\begin{align*}
    \|h\|_\infty &= \sum_{k \geq 1} |\langle h,\varphi_k\rangle_{L_2}| \cdot \|\varphi_k\|_\infty \\
    & \leq \sqrt{2/\pi} \cdot \sqrt{L} \left(\sum_{j, j' \geq 0, \, (j,j') \not = (0,0)} (j^2 + (j')^2)^{-s} \right)^{1/2} \leq C_0,
\end{align*}
where $C_0=C_0(s,L)$ only depends on $s>1$ and $L>0$.

We show analogously that $h$ is Lipschitz in its first argument whenever $f$ and $\tilde f$ belong to $W(s,L,\rho^*)$ with $s>2$. Indeed,
\begin{align*}
    \|\frac{\partial}{\partial t} h(t,x)\|_\infty &= \sum_{j,j' \geq 0} |\langle h,\varphi_{j,j'}\rangle_{L_2}| \cdot 2\pi j \|\varphi_{j,j'}\|_\infty \\
    &= 2\pi \sqrt{2/\pi} \sum_{j,j' \geq 0, \, (j,j') \not = (0,0)} |\langle h,\varphi_{j,j'}\rangle_{L_2}| \cdot ( j^2 +(j')^2)^{1/2} \\
    & \leq 2\sqrt{2\pi} \cdot \sqrt{L} \left(\sum_{j, j' \geq 0, \, (j,j') \not = (0,0)} (j^2 + (j')^2)^{-s+1}) \right)^{1/2} \leq C_1,
\end{align*}
where $C_1 = C_1(s,L)$ only depend on $s>2$ and $L>0$. 

As the map $f \mapsto \theta(f)$ is linear we have, for $h=f-\tilde{f}$, that
\begin{align*}
	& \|\theta(f) - \theta(\tilde{f})\|_F^2  \, = \|\theta(h)\|_F^2 \\ & \, = \, 2 \sum_{j=0}^{n-1} \sum_{j'=0}^{j-1} \Big| \int_{-\pi}^\pi \exp(i(j'-j)x) h(j'/n,x) dx\Big|^2 \, + \, \sum_{j=0}^{n-1} \Big|\int_{-\pi}^\pi h(j/n,x) dx\Big|^2 \\
	& \, = \, 4\pi \sum_{j'=0}^{n-1} \sum_{j=j'+1}^{n-1} \Big| (2\pi)^{-1/2} \int_{-\pi}^\pi \exp(-ijx) \exp(ij'x) h(j'/n,x) dx\Big|^2  \\ & \hspace{9cm} +  \sum_{j=0}^{n-1} \Big|\int_{-\pi}^\pi h(j/n,x) dx\Big|^2 \\ 
	& \, \leq \, 6\pi n \cdot \frac 1n \sum_{j=0}^{n-1} \int_{-\pi}^\pi |h(j/n,x)|^2 dx \leq 12 \pi^2 n \|h\|_\infty^2 ,
\end{align*}
by Parseval's identity and the Cauchy-Schwarz inequality. 

We show that the Riemann sum in the previous line can be fast approximated by the integral $\|h\|^2_{L_2}$. Indeed,
\begin{align*}
    \frac 1n \sum_{j=0}^{n-1} \int_{-\pi}^\pi & |h(j/n,x)|^2 dx - \|h\|^2_{L_2}
    =  \sum_{j=0}^{n-1} \int_{j/n}^{(j+1)/n} \int_{-\pi}^\pi \left( |h(j/n,x)|^2  - h(t,x)^2 \right) dx\, dt\\
    & \leq 2 \|h\|_\infty  \sum_{j=0}^{n-1} \int_{j/n}^{(j+1)/n} \int_{-\pi}^\pi \left| h(j/n,x)  - h(t,x) \right| dx\, dt\\
   & \leq 2  \|h\|_\infty C_1  \sum_{j=0}^{n-1} \int_{j/n}^{(j+1)/n} \int_{-\pi}^\pi \left| j/n - t \right| dx\, dt \leq \frac{4 \pi C_1  \|h\|_\infty}{n},
\end{align*}
which yields the desired result. \hfill $\square$ \\[0.5cm]

\noindent {\it Proof of Lemma~\ref{lem:pilot}}: It is sufficient to show that both
$$
\mathbb{E} \left[ \| \hat \alpha - a\|^2 \right] \leq \mathcal{O}(1) \cdot K_n
\text{ and }
 \|a- \alpha_{\theta(f)} \|^2  \leq \mathcal{O}(1) \cdot K_n.
$$
First, note that
$$
 \| a- \alpha_{\theta(f)} \|^2 = \sum_{j=1}^{K_n} a_j^2 \left(1 - \frac{\|M_j\|_F}{\sqrt{2 \pi n}} \right)^2 
 \leq \frac 2{\pi} \sum_{j=1}^{K_n} \left(1 - \sqrt{ 1 - \frac{j}{n }} \right)^2 
 \leq \frac 2{\pi} \sum_{j=1}^{K_n} \frac{\kappa_2^2}{ n^2} \leq \mathcal{O}\left( \frac{K_n^2}{n^2}\right), 
$$
which is $\mathcal{O}(1) \cdot K_n$ as soon as $K_n \leq \mathcal{O}(1) \cdot n^2$. 
Then, let us see that
$$
\mathbb{E} \left[ \| \hat  \alpha - a\|^2 \right] = 2 \pi n \cdot \sum_{j=1}^{K_n} \mathbb{E} \,\left[ \left| \langle \hat f, \varphi_j \rangle_{L_2} - \langle f, \varphi_j \rangle_{L_2} \right|^2 \right]. 
$$
We have, for an arbitrary constant $c>0$:
\begin{align}
   & \mathbb{E} \,\left[ \left| \langle \hat f, \varphi_j \rangle_{L_2} - \langle f, \varphi_j \rangle_{L_2} \right|^2 \right]
   =  \mathbb{E} \left[ \left|\int \left( e^{\widetilde{\log } f } -  f\right) \varphi_j \right|^2 \right] \nonumber \\
   & \leq \mathbb{E} \left[ \left|\int \left( e^{\widetilde{\log } f - \log f} -  1 \right) \cdot f \varphi_j \right|^2 \right] \nonumber \\
   & \leq \mathbb{E} \left| \sum_{k \geq 1} \frac {c^{k/2} \cdot c^{-k/2}}{k!} \int \left( \widetilde{\log } f - \log f\right)^k \cdot f \varphi_j \right|^2 \nonumber \\
   & \leq \sum_{k \geq 1}  \frac {c^k}{k!} \cdot \sum_{k \geq 1}  \frac {c^{-k} }{k!} \mathbb{E} \left[ \langle \left( \widetilde{\log } f - \log f\right)^k  ,  f \varphi_j\rangle^2 \right] \label{eq1}
\end{align}
where we used Cauchy-Schwarz. We treat separately the term for $k=1$ and the sum over $k \geq 2$.   We decompose  
$$
\frac{\tilde \alpha_\ell}{\sqrt{2\pi n} } = \langle \varphi_\ell , \log f \rangle_{L_2} + a_n  \cdot \langle \varphi_\ell , dW \rangle_{L_2}.
$$
where we denote $a_n := 2\sqrt{\pi/n}$. For $k=1$:
\begin{align*}
     \langle \left( \widetilde{\log } f - \log f\right)  ,  f \varphi_j\rangle
     & =  \sum_{\ell > J_n} \langle \varphi_\ell , \log f \rangle_{L_2} \cdot  \int \varphi_\ell \varphi_j f + a_n \sum_{\ell = 1}^{J_n}  \langle \varphi_\ell ; dW \rangle_{L_2}  \cdot \int \varphi_\ell \varphi_j f \, .
\end{align*}
Use that $|\int \varphi_k \varphi_j f| \leq 4 /\rho^*$ and that the random variables $ \langle \varphi_k;dW \rangle_{L_2} $ for $k=1,\ldots,K_n$ are centered and independent to get
\begin{align}
     \mathbb{E} \left[ \langle \left( \widetilde{\log } f - \log f\right)  ,  f \varphi_j\rangle^2 \right]
     & \leq \left( \frac 4{\rho^*}  \sum_{\ell > J_n} \langle \varphi_\ell , \log f\rangle_{L_2} \right)^2 +  a_n^2 \sum_{\ell = 1}^{K_n} ( \int \varphi_\ell \varphi_j f )^2 \nonumber \\
     & \leq \left( \frac 4{\rho^*}  \sum_{\ell > J_n} \left| \langle \varphi_\ell , \log f \rangle_{L_2}\right| \right)^2 +  a_n^2  \int ( \varphi_j f )^2 \nonumber \\
     & \leq \frac{16}{ (\rho^*)^2}\, B(J_n)^2 +  a_n^2  \frac{4}{(\rho^*)^2} \label{eq2}
\end{align}
where we denote $B(J_n):=  \sum_{\ell \geq J_n} \left| \langle \varphi_\ell , \log f \rangle_{L_2}\right| $. For $k \geq 2$:
\begin{align}
     & \mathbb{E} \left[ \langle \left( \widetilde{\log } f - \log f\right) ^k ,  f \varphi_j\rangle^2 \right]
      \leq \int (\varphi_j f)^2 \cdot \mathbb{E} \left[ \int \left( \widetilde{\log } f - \log f\right)^{2k} \right] \nonumber \\
     & \leq \frac 4{(\rho^*)^2} \cdot \mathbb{E} \left[ \int \left( \sum_{\ell > J_n} \varphi_\ell \langle  \varphi_\ell , \log f \rangle_{L_2} + a_n \sum_{\ell =1}^{J_n} \varphi_\ell \langle \varphi_\ell ; dW \rangle_{L_2}  \right)^{2k} \right] \nonumber \\
     & \leq \frac {4 \cdot 2^{2k-1} }{(\rho^*)^2} \left[ \int \left(\sqrt{\frac{2}{\pi}} B(J_n) \right)^{2k} + a_n^{2k} \int \mathbb{E} \left[\left( \sum_{\ell =1}^{J_n} \varphi_\ell \langle \varphi_\ell ; dW \rangle_{L_2} 
 \right)^{2k} \right] \right] \nonumber \\
      & \leq \frac {2 \cdot 2^{2k} }{(\rho^*)^2} \left[ 2 \pi (\frac{2}{\pi} ) ^k  B(J_n) ^{2k} + a_n^{2k} \int \mathbb{E} \left[\left( \sum_{\ell =1}^{J_n} \varphi_\ell \langle \varphi_\ell ; dW \rangle_{L_2} 
 \right)^{2k} \right] \right] \, . \label{eq3}
\end{align}
Note that $\xi(t,x) :=\sum_{\ell =1}^{J_n} \varphi_\ell(t,x) \langle \varphi_\ell ; dW \rangle_{L_2} $ is a Gaussian random variable with mean  0 and variance $ \sum_{\ell =1}^{J_n} \varphi_\ell^2(t,x)  . $ We deduce that 
$$
\mathbb{E} \xi^{2k} = \frac{(2k-1)! }{ 2^{k-1} (k-1)!} \cdot (\sum_{\ell =1}^{J_n} \varphi_\ell^2(t,x))^k 
\leq \frac{(2k-1)! }{ 2^{k-1} (k-1)!} \cdot \left(\frac 2{\pi} J_n \right)^k\, .
$$
Plug this bound into \eqref{eq3}, then \eqref{eq3} and \eqref{eq2} into \eqref{eq1} to get:
\begin{align*}
    & \mathbb{E} \,\left[ \left| \langle \hat f, \varphi_j \rangle_{L_2} - \langle f, \varphi_j \rangle_{L_2} \right|^2 \right]
    \leq \frac{ 16 e^c}{c (\rho^*)^2} B(J_n)^2 + a_n^2 \frac{4 e^c}{(\rho^*)^2}  \\
    & + \frac{4 \pi e^c}{(\rho^*)^2} \sum_{k \geq 2} \frac 1{k!} \left( \frac{8 B(J_n)^2}{c \pi}\right)^k + \frac{8 \pi e^c}{(\rho^*)^2} \sum_{k \geq 2} {{2k-1}\choose{k}} \left( \frac{4 a_n^2 J_n}{c \pi } \right)^k\\
    & \leq \mathcal{O}(1) \left( B(J_n)^2 + a_n^2 + B(J_n)^4 + a_n^4 J_n^2 \right) \, ,
\end{align*}
for $c>0$ large enough and using that ${{2k-1}\choose {k}} \leq 5^k$ for $k\geq 2$ large enough. In order to finish the proof let us show an upper bound on $B(J_n)$ that allows us to choose $J_n$. 

The spectral density $f$ belongs to $W(s,L,\rho^*)$ and therefore it has weak partial derivatives $\partial^{j+j'}/(\partial^jt \partial^{j'}x) f(t,x)$ for all non-negative $j,\, j'$ such that $j+j' \leq \lfloor s \rfloor$ and they are all periodic and uniformly bounded in $L_2$ norm. It is easy to see that $\log f$ is a bounded function and has weak partial derivatives for all positive $j,\, j'$ such that $j+j' \leq s^*$ which are all uniformly bounded in $L_2$ norm, for some $s^* \leq s+1$. Thus,
$$
 \left| \langle \varphi_{k,k'} , \log f \rangle_{L_2}\right| \leq \int \frac{1}{(2\pi k)^j \cdot (k')^{j'}} |\varphi_{k,k'}| \left| \frac{\partial^{j+j'}}{\partial^j t \, \partial^{j'}x} \log f(t,x)\right| dt dx \leq \frac{C}{k^j (k')^{j'}},
$$
where the positive constant $C$ depends only on $s,\, L$ and $\rho^*$. Thus, let us take 
\begin{align*}
    B(J_n) & = \sum_{k,k'}  \, {\bf 1}\{k>\kappa_1\mbox{ or }k'>\kappa_2\} \, \left| \langle \varphi_{k,k'} , \log f \rangle_{L_2}\right| \\
    & \leq  \sum_{k,k'}  \, {\bf 1}\{k>\kappa_1\mbox{ or }k'>\kappa_2\} \frac{C}{k^j (k')^{j'}} \\
    &\leq \mathcal{O}(1) \cdot \kappa_1^{-j+1} \vee \kappa_2^{-j'+1} 
    = \mathcal{O}(1) \cdot J_n^{-\lfloor {s^*}/{2} \rfloor +1}\, ,
\end{align*}
where we choose $\kappa_1 = \kappa_2 = J_n^{1/2}$ and $j=j'= \lfloor {s^*}/{2} \rfloor$.
Finally,
$$
 \mathbb{E} \,\left[ \left| \langle \hat f, \varphi_j \rangle_{L_2} - \langle f, \varphi_j \rangle_{L_2} \right|^2 \right]  \leq \mathcal{O}(1) \left( J_n^{-2 \lfloor {s^*}/{2} \rfloor +2} + a_n^2 + a_n^2 (a_n J_n)^2 \right) \, 
$$
which becomes $\mathcal{O}( 1/n)$ by choosing $J_n = n^{1/2}$ and for $ \lfloor {s^*}/{2} \rfloor -1 \geq 1$, that is $s^* \geq 4$.
\hfill $\Box$ \\
[1cm]

\noindent {\it Proof of Lemma \ref{L:f12}}: The function $\hat{f}^{-1/2}$ is real-valued, bounded, measurable and symmetric in the second component; thus, it is located in the corresponding Hilbert-subspace of $L_2([0,1]\times [-\pi,\pi])$ so that it can be expanded with respect to the orthonormal basis $\{\varphi_j\}_j$ such that
$$ \hat{f}^{-1/2} \, = \, \sum_{j=1}^\infty \langle \hat{f}^{-1/2} , \varphi_j\rangle_{L_2} \cdot \varphi_j\,. $$
By $\partial_{\ell,\ell'}$ we denote the partial differential operator $\partial^{\ell+\ell'}/(\partial t^\ell \, \partial x^{\ell'})$. Note that
\begin{align*}
\big\|\partial_{\ell,\ell'} f_n\big\|_\infty & \, \leq \, \sum_{j,j'=0}^{\kappa_1,\kappa_2} |\tilde{f}^+_{j,j'}| \cdot  \|\partial_{\ell,\ell'}\varphi^+_{j,j'}\|_\infty \, + \, \sum_{j=1,j'=0}^{\kappa_1,\kappa_2} |\tilde{f}^-_{j,j'}| \cdot  \|\partial_{\ell,\ell'}\varphi^-_{j,j'}\|_\infty \\
& \, \leq \, \sqrt{2/\pi} \sum_{j,j'=0}^{\kappa_1,\kappa_2} \big(|\tilde{f}^+_{j,j'}| \, + \, |\tilde{f}^-_{j,j'}|\big) \cdot (2\pi j)^\ell (j')^{\ell'} \\
& \, \leq \, (2/\sqrt{\pi})\cdot (2\pi)^{\ell } \cdot \sqrt{L} \cdot\Big(\sum_{j,j'} {\bf 1}\{j^2+(j')^2>0\} \cdot \big(j^2+(j')^2\big)^{\ell+\ell'-s}\Big)^{1/2}\,, 
\end{align*}
whenever $1\leq \ell+\ell'\leq s^* < s- 1$. Under this condition $\big\|\partial_{\ell,\ell'} f_n\big\|_\infty$ admits an upper bound that only depends on $L$ and $s$. On the other hand, for $\ell=\ell'=0$, an upper bound has already been provided by $2/\rho^*$. Furthermore,
\begin{align*}
\big\|\partial_{\ell,\ell'} (\hat{f} - f_n)\big\|_\infty & \, \leq \, \frac1{\sqrt{2\pi n}} \, \sum_{k=1}^{K_n} |\hat \alpha_k - \alpha_{\theta(f),k}| \cdot \|\partial_{\ell,\ell'} \varphi_k\|_\infty \\
& \, \leq \, \frac{\sqrt{K_n}}{\pi \sqrt{ n}} \, \Big(\sum_{k=1}^{K_n} |  \hat \alpha_k - \alpha_{\theta(f),k} |^2\Big)^{1/2}\cdot (2\pi \kappa_1)^\ell\cdot \kappa_2^{\ell'} \\
& \, \leq \, \frac{\sqrt{K_n}\gamma_n}{\pi \sqrt{n}} \cdot (2\pi \kappa_1)^\ell\cdot \kappa_2^{\ell'} \, \asymp \, K_n^{(\ell+\ell'+1)/2} \cdot \gamma_n \cdot n^{-1/2}\,, 
\end{align*}
by (\ref{eq:gamman}). Under the conditions of the lemma we have established an upper bound on $\|\partial_{\ell,\ell'} \hat{f}\|_\infty$ only depending on $s^*$, $s$ and $L$ as well when $\ell+\ell'\leq s^*$. Now we introduce the differential operators $\partial_b$, $b\in \{1,\ldots,\ell+\ell'\}$, where all of these operators are either $\partial/(\partial t)$ or $\partial/(\partial x)$; and the notation $\partial_B := \prod_{b\in B} \partial_b$. By the multivariate version of Fa\`{a} di Bruno's formula we deduce that
\begin{align*}
\partial_{\ell,\ell'} \hat{f}^{-1/2} & \, = \, 2\, \sum_P (-1/4)^{\# P} \frac{(2\# P - 1)!}{(\# P -1)!} \, \hat{f}^{-1/2-\#P} \cdot \prod_{B \in P} \partial_B \hat{f}\,,
\end{align*}
where the sum is taken over all partitions $P$ of the set $\{1,\ldots,\ell+\ell'\}$. By the derived upper bounds on the supremum norm of the partial derivatives of $\hat{f}$ and since $\hat{f} \geq \rho^*/2>0$ we attain a uniform upper bound on $\|\partial_{\ell,\ell'} \hat{f}^{-1/2}\|_\infty$ for $\ell+\ell'\leq s^*$ which only depends on $s$, $s^*$, $L$ and $\rho^*$. Moreover the periodicity of $\hat{f}$, i.e. the property that $\hat{f}(t+2\pi k,x+k') = \hat{f}(t,x)$ for all $t,x\in \mathbb{R}$ and integer $k$ and $k'$, extends to the function $\hat{f}^{-1/2}$ and its partial derivatives. Therefore, using integration by parts with vanishing boundary terms,
\begin{align*} 
\langle \hat{f}^{-1/2} , \varphi_j \rangle_{L_2} & \, = \, (-1)^{\ell+\ell'}\cdot \big\langle \partial_{\ell,\ell'} \hat{f}^{-1/2} , \partial_{\ell,\ell'}^{-1} \varphi_j \big\rangle_{L_2}\,,
\end{align*}
where $\|\partial_{\ell,\ell'}^{-1} \varphi_{k,k'}^\pm\|_\infty \leq \sqrt{2/\pi} \cdot (2\pi k)^{-\ell}\cdot (k')^{-\ell'}$ while we stipulate that $\ell=0$ if $k=0$ 
and $\ell'=0$ in the case of $k'=0$. For all $(k,k')\neq (0,0)$ we choose $\ell$, $\ell'$ such that $\ell+\ell'=s^*$. Then,
\begin{align*} \big\|\widetilde{f^{-1/2}} - & \hat{f}^{-1/2}\big\|_\infty \\ &  \leq  \sqrt{2/\pi} \, \sum_{k,k'} {\bf 1}\{k>\kappa_1\mbox{ or }k'>\kappa_2\}  \big(\big|\langle \hat{f}^{-1/2} , \varphi_{k,k'}^+\rangle_{L_2}\big| + \big|\langle \hat{f}^{-1/2} , \varphi_{k,k'}^-\rangle_{L_2}\big|\big) \\ &   \asymp \kappa_1^{1-\ell} \cdot \kappa_2^{1-\ell'} \, \asymp \, K_n^{1-s^*/2}\,, \end{align*}
completes the proof of the lemma. \hfill $\square$ \\[5mm]

\noindent {\it Proof of Lemma \ref{L:almosthom}}: We have
  \begin{align*}
  A & \, = \, \sum_{j=-\kappa_1}^{\kappa_1} \sum_{j'=-\kappa_2}^{\kappa_2} a_{j,j'} \, \check{M}_{0,j,j'}\,, \\
  B & \, = \, \sum_{j=-\kappa_1}^{\kappa_1} \sum_{j'=-\kappa_2}^{\kappa_2} b_{j,j'} \, \check{M}_{0,j,j'}\,,
  \end{align*}
  for some coefficients $a_{j,j'}$ and $b_{j,j'}$. The orthogonality of the $\check{M}_{j,j'}$ yields that
  \begin{align*}
  \|A\|_F^2 & \, = \, n \cdot \sum_{j=-\kappa_1}^{\kappa_1} \sum_{j'=-\kappa_2}^{\kappa_2} |a_{j,j'}|^2\,, \\
  \|B\|_F^2 & \, = \, n \cdot \sum_{j=-\kappa_1}^{\kappa_1} \sum_{j'=-\kappa_2}^{\kappa_2} |b_{j,j'}|^2\,.
  \end{align*}
  Using (\ref{eq:almosthom}) we deduce that
  \begin{align*}
  \Psi(AB) & \, = \, \sum_{j,j'} \sum_{k,k'} a_{j,j'} b_{k,k'} \, \Psi(\check{M}_{0,j,j'} \check{M}_{0,k,k'}) \\
  & \, = \, \sum_{j,j'} \sum_{k,k'} a_{j,j'} b_{k,k'} \, \exp(2\pi i j' k/n) \, \Psi(\check{M}_{0,j+k,j'+k'}) \\
  & \, = \, \sum_{j,j'} \sum_{k,k'} a_{j,j'} b_{k,k'} \, \exp(2\pi i j' k/n) \, \tilde{\varphi}_{j+k,j'+k'}\,.
  \end{align*}
  On the other hand,
  \begin{align*}
  \Psi(A) \Psi(B) & \, = \, \sum_{j,j'} \sum_{k,k'} a_{j,j'} b_{k,k'} \, \Psi(\check{M}_{0,j,j'}) \Psi(\check{M}_{0,k,k'}) \\
  & \, = \, \sum_{j,j'} \sum_{k,k'} a_{j,j'} b_{k,k'} \, \tilde{\varphi}_{j+k,j'+k'}\,.
  \end{align*}
  Therefore,
  \begin{align*}
  & \, = \, \sum_{\ell,\ell'} \Big|\sum_{j,j'} a_{j,j'} b_{\ell-j,\ell'-j'} \big(\exp(2\pi i j'(\ell-j)/n\big)-1\big)\Big|^2 \big\|\tilde{\varphi}_{\ell,\ell'}\|_{L_2^n}^2 \\
  & \, \leq \, n \cdot \sum_{j,j'} |a_{j,j'}|^2 \sum_{\ell,\ell'} |b_{\ell-j,\ell'-j'}|^2 \big|\exp(2\pi i j'(\ell-j)/n)-1\big|^2 \\
  & \, \leq \, n \cdot \sum_{j,j'} |a_{j,j'}|^2 \sum_{k,k'} |b_{k,k'}|^2 \big(2 - 2\cos(2\pi j'k/n)\big) \\
  & \, \leq \, n \cdot \sum_{j,j'} |a_{j,j'}|^2 \sum_{k,k'} |b_{k,k'}|^2 \cdot 4\pi^2 (j'k)^2/n^2  \\
  & \, \leq \, 4\pi^2 \, \|A\|_F^2 \|B\|_F^2 \cdot \kappa_1^2 \kappa_2^2 / n^3\,,
  \end{align*}
  by the orthogonality of the $\tilde{\varphi}_{j,j'}$ and the Cauchy-Schwarz inequality where the sums may be taken over all integers when the coefficients $a_{j,j'}$ and $b_{j,j'}$ are put equal to $0$ for $|j|>\kappa_1$ or $|j'|>\kappa_2$. \hfill $\square$ \\[0.5cm]

\newpage


\bibliographystyle{imsart-nameyear} 
\bibliography{sn-LSAE}       

\end{document}